\documentclass[11pt]{article}

\newcommand{\documentdate}{3 November 2018}

\usepackage{a4wide,latexsym,varioref,amsmath}

\newcommand{\numsection}[1]{\section{#1}\setcounter{equation}{0}}

\renewcommand{\theequation}{\arabic{section}.\arabic{equation}}

\newcommand{\beqn}[1]{\begin{equation}\label{#1}}
\newcommand{\eeqn}{\end{equation}}
\newcommand{\req}[1]{(\ref{#1})}
\newcommand{\tim}[1]{\;\; \mbox{#1} \;\;}
\newcommand{\flow}{f_{\rm low}}
\DeclareMathOperator*\globmin{globmin}
\DeclareMathOperator*\globmax{globmax}

\newcommand{\bpr}{{\bf Proof.} \hspace{1.5mm}}
\newcommand{\epr}{\hfill $\Box$ \vspace*{1em}}
\newcommand{\eqdef}{\stackrel{\rm def}{=}}
\newcommand{\bigfrac}[2]{\frac{\displaystyle #1}{\displaystyle #2}}

\newcommand{\bigsum}{\displaystyle \sum}
\newcommand{\bigint}{\displaystyle \int}

\newcommand{\bigmax}{\displaystyle \max}
\newcommand{\calS}{{\cal S}}
\newcommand{\calF}{{\cal F}}
\newcommand{\calU}{{\cal U}}
\newcommand{\ms}{\;\;\;\;}
\newcommand{\sfrac}[2]{{\scriptstyle \frac{#1}{#2}}}
\newcommand{\half}{\sfrac{1}{2}}
\newcommand{\threequarters}{\sfrac{3}{4}}
\newcommand{\ii}[1]{\{1, \ldots, #1 \}}
\newcommand{\iibe}[2]{\{ #1, \ldots, #2 \}}
\newcommand{\mystack}[2]{_{\stackrel{\scriptstyle #1}{\scriptstyle #2}}}
\newcommand{\matr}[2]{\left( \begin{array}{#1} #2 \end{array} \right) }
\newcommand{\cvect}[1]{\left( \begin{array}{c} #1 \end{array} \right) }
\renewcommand{\Re}{\hbox{I\hskip -2pt R}}

\newcommand{\Na}{\hbox{I\hskip -1.8pt N}}
\newcommand{\proof}[1]{
\begin{list}{}{
\setlength{\topsep}{0.0pt}
\setlength{\partopsep}{0.0pt}
\setlength{\leftmargin}{0.025\textwidth}
\setlength{\rightmargin}{0.5\leftmargin}
\setlength{\labelwidth}{0.5\leftmargin}
\setlength{\labelsep}{0.25\leftmargin}}
\item \bpr #1 \epr \noindent
\end{list}}
\newtheorem{theorem}{Theorem}[section]
\newtheorem{lemma}[theorem]{Lemma}

\newcommand{\llem}[2]{\vspace{\baselineskip} 
\noindent\framebox[\textwidth]{\parbox{0.95\textwidth}{
\begin{lemma} \label{#1} \rm #2 \end{lemma} } } \vspace{\baselineskip} }
\newcommand{\lthm}[2]{\vspace{\baselineskip} 
\noindent\framebox[\textwidth]{\parbox{0.95\textwidth}{
\begin{theorem} \label{#1} \rm #2 \end{theorem} } } \vspace{\baselineskip} }
\newcounter{algo}[section]
\renewcommand{\thealgo}{\thesection.\arabic{algo}}
\newcommand{\algo}[3]{\refstepcounter{algo}
\begin{center}\begin{figure}[htbp]
\framebox[\textwidth]{
\parbox{0.95\textwidth} {\vspace{\topsep}
{\bf Algorithm \thealgo : #2}\label{#1}\\
\vspace*{-\topsep} \mbox{ }\\
{#3} \vspace{\topsep} }}
\end{figure}\end{center}}

\title{Sharp worst-case evaluation complexity bounds for arbitrary-order 
  nonconvex optimization with inexpensive constraints}

\author{
C. Cartis\thanks{Mathematical Institute,
   Oxford University,
   Oxford OX2 6GG, England.  Email: coralia.cartis@maths.ox.ac.uk},
N. I. M. Gould\thanks{Computational Mathematics Group,
   STFC-Rutherford Appleton Laboratory,
   Chilton OX11 0QX, England. Email:  nick.gould@stfc.ac.uk .
   The work of this author was supported by EPSRC grant EP/M025179/1}
~and~Ph. L. Toint\thanks{Namur Center for Complex Systems (naXys),
   University of Namur, 61, rue de Bruxelles, B-5000 Namur, Belgium.
   Email: philippe.toint@unamur.be} 
}

\date{\documentdate}

\begin{document}

\maketitle

\begin{abstract}
We provide sharp worst-case evaluation complexity bounds for nonconvex
minimization problems with general inexpensive constraints, i.e.\ problems
where the cost of evaluating/enforcing of the (possibly nonconvex or even
disconnected) constraints, if any, is negligible compared to that of
evaluating the objective function. These bounds unify, extend or improve all
known upper and lower complexity bounds for unconstrained and
convexly-constrained problems. It is shown that, given an accuracy level
$\epsilon$, a degree of highest available Lipschitz continuous derivatives $p$
and a desired optimality order $q$ between one and $p$, a conceptual
regularization algorithm requires no more than
$O(\epsilon^{-\frac{p+1}{p-q+1}})$ evaluations of the objective function and
its derivatives to compute a suitably approximate $q$-th order minimizer. With
an appropriate choice of the regularization, a similar result also holds if
the $p$-th derivative is merely H\"{o}lder rather than Lipschitz
continuous. We provide an example that shows that the above complexity bound
is sharp for unconstrained and a wide class of constrained problems; we also
give reasons for the optimality of regularization methods from a worst-case
complexity point of view, within a large class of algorithms that use the same
derivative information.
\end{abstract}
 
\numsection{Introduction}

Since the seminal paper by Vavasis \cite{Vava93} on the complexity of finding
first-order critical points in unconstrained nonlinear optimization was
published 25 years ago, the question of the optimal worst-case complexity of
optimization methods has been of interest to mathematicians and also, because
of its strong connection with deep learning, to computer scientists. Of late,
there has been a growing interest in this research field, both for convex and
nonconvex problems. This paper focusses on the latter class and follows a now
subtantial\footnote{See \cite{CartGoulToin18a} for a more complete list of
references.} trend of research where bounds on the worst-case evaluation
complexity (or oracle complexity) of obtaining first- and (more rarely)
second-order-necessary minimizers\footnote{That is points satisfying the
first- or second-order necessary optimality conditions for minimization.}
for nonlinear nonconvex unconstrained optimization problems
\cite{Vava93,Nest04,GratSartToin08,NestPoly06,CartGoulToin11d}.  These papers
all provide \emph{upper} evaluation complexity bounds: they show that, to
obtain an $\epsilon$-approximate first-order-necessary minimizer (for
unconstrained problem, this is a point at which the gradient of the objective
function is less than $\epsilon$ in norm), \emph{at most} $O(\epsilon^{-2})$
evaluations of the objective function\footnote{And its available derivatives.}
are needed if a model involving first derivatives is used, and \emph{at most}
$O(\epsilon^{-3/2})$ evaluations are needed if using second derivatives is
permitted. This result was extended to convexly-constrained problems in
\cite{CartGoulToin12b}. A broader framework allowing the use of Taylor series
of degree $p$ was more recently proposed in \cite{BirgGardMartSantToin17}, in
which case the worst-case evaluation complexity bound for
$\epsilon$-first-order-necessary unconstrained minimizer is shown to be
$O(\epsilon^{-\frac{p+1}{p}})$, thereby generalizing the previous results for
this case.  Complexity for obtaining $\epsilon$-approximate
second-order-necessary unconstrained minimizers was considered in
\cite{NestPoly06,CartGoulToin11d}, where a bound of $O(\epsilon^{-3})$
evaluations was proved to obtain an $\epsilon$-second-order-necessary
minimizer using a Taylor's model of degree two, and a bound of
$O(\epsilon^{-\frac{p+1}{p-1}})$ evaluations was shown in
\cite{CartGoulToin17e} for the case where a Taylor model of degree $p$ is
used. Defining $q$-th-order-necessary minimizers for $q > 2$ was considered in
\cite{CartGoulToin17c}, where the difficulty of stating and verifying
necessary optimality was discussed. In particular, it was concluded in this
latter reference that defining and computing $\epsilon$-approximate
$q$-th-order-necessary minimizers for $q>2$ is likely to remain elusive,
essentially because of the nonlinearity and lack of continuity of the kernels
of the derivatives involved. A more general Taylor-based definition of
optimality was introduced instead, which allowed to show an upper bound of
$O(\epsilon^{-(q+1)})$ on evaluation complexity for convexly-constrained
problems, in particular improving on the bound of $O(\epsilon^{-9/2})$ stated
in \cite{AnanGe16} for the case $p=q=3$.

The unconstrained and convexly-constrained cases where the assumption of
Lipschitz continuity is replaced by the weaker $\beta$-H\"{o}lder continuity
($\beta \in (0,1]$) have also been studied for $q=1$ in
\cite{NestGrap16,CartGoulToin16,CartGoulToin17d}.  These references
show that \emph{at most} $O(\epsilon^{-\frac{p+\beta}{p-1+\beta}})$
evaluations are needed for obtaining an $\epsilon$-first-order-necessary
minimizer.

While upper complexity bounds are important as they provide a handle on the
intrinsic difficulty of the considered problem, they do so at the condition of
not being overly pessimistic. To address this last point, \emph{lower} bounds
on the evaluation complexity of unconstrained nonconvex optimization problems
and methods were derived in \cite{CartGoulToin10a,Nest04}  and
\cite{CartGoulToin18a}, where it was shown that the known upper complexity 
bounds are sharp (irrespective of problem's dimension) for most known methods
using Taylor's models of degree one or two. That is to say that there are
examples for which the complexity order predicted by the upper bound is
actually achieved. More recently, Carmon \emph{et al.}
\cite{CarmDuchHindSidf17a} provided an elaborate construction showing that
\emph{at least} a multiple of $\epsilon^{-\frac{p+1}{p}}$ function evaluations
may be needed to obtain an $\epsilon$-first-order-necessary unconstrained
minimizer where derivatives of order at most $p$ are used. This result, which
matches in order the upper bound of \cite{BirgGardMartSantToin17}, covers a
very wide class of potential optimization methods\footnote{In particular, it
covers randomized   methods, which we do not consider in this paper.} but has
the drawback of being only valid for problems whose dimension essentially
exceeds the number of iterations needed, which can be very large and quickly
grows when $\epsilon$ tends to zero.

{\bf Contributions.}  The present paper aims at unifying and generalizing all
the above results in a single framework, providing, for problems with
inexpensive or no constraints, provably optimal evaluation complexity bounds
for arbitrary optimality order, all relevant model degrees and levels of
smoothness of the objective function. By ``inexpensive constraints'', we mean
general set constraints whose enforcement and evaluation\footnote{Constraint's
  values and that of their derivatives, if relevant.} cost is negligible
compared to the cost of evaluating the objective function. As a consequence,
the evaluation complexity for such problems is meaningfully captured by focusing
of the number of evaluations of this latter function. This class of
minimization problems contains important cases such as bound-constrained
problems and convexly-constrained problems (when the projection onto the
feasible set is inexpensive), but also allows possibly nonconvex or even
disconnected feasible sets.

In order to achieve these objectives, we first revisit the Taylor-based
optimality measure of \cite{CartGoulToin17c} and define
($\epsilon$,$\delta$)-$q$-th-order-necessary minimizers, a notion extending
the standard $\epsilon$-first- and $\epsilon$-second-order cases to arbitrary
orders.  We then present a conceptual regularization algorithm using degree
$p$ models and show that this algorithm requires at most
$O(\epsilon^{-\frac{p+\beta}{p-q+\beta}})$ evaluations of $f$ and its
derivatives to find such an ($\epsilon$,$\delta$)-$q$-th-order-necessary
minimizer when the $p$-th derivative of $f$ is assumed to be
$\beta$-H\"{o}lder continuous. (If the $p$-th derivative is assumed to be
Lispchitz continuous, the bound becomes $O(\epsilon^{-\frac{p+1}{p-q+1}})$.)
This bound matches the best known lower bounds for first- and second-order,
and improves on the bound in $O(\epsilon^{-(q+1)})$ given by
\cite{CartGoulToin17c}. We then show that this bound is sharp in order for
unconstrained problems with Lipschitz continuous $p$-th derivative by
completing and extending the result of \cite{CarmDuchHindSidf17a} in two
ways. The first is to show that the lower worst-case bound of order
$\epsilon^{-\frac{p+1}{p}}$ evaluations for obtaining a first-order-necessary
minimizer using at most $p$ derivatives is also valid for problems of every
dimension, and the second is to show that this bound can be generalized to a
multiple of $\epsilon^{-\frac{p+1}{p-q+1}}$ for obtaining a
$q$-th-order-necessary minimizer of any order $q$. In particular, this result
matches in order the upper bound obtained in the first part of the paper and
subsumes or improves known lower bounds for first- and second-order-necessary
minimizers.  While our lower bounds are derived for regularization algorithms
applied to unconstrained problems, we also indicate that they may be extended
to a much wider class of minimization methods and to a significant class of
constrained problems.

The paper is organized as follows. Section~\ref{arp-s} introduces the
(possibly constrained) minimization problem of interest and the concept of
($\epsilon$,$\delta$)-approximate $q$-th-order-necessary minimizers. It also
presents a variant of the Adaptive Regularization algorithm using degree $p$
Taylor's models (AR$p$) whose purpose is to find such minimizers.
Section~\ref{upperbound-s} then provides an upper bound on the evaluation
complexity for the AR$p$ algorithm to achieve this
task. Section~\ref{2ndorder-s} then discusses
specialization of this result to the case where $\epsilon$-approximate
second-order-necessary minimizers are sought.  The complexity upper bound of
Section~\ref{upperbound-s} is then proved to be sharp in
Section~\ref{lowerbound-s} for the Lipschitz-continuous cases where the feasible set contains a
ray. Some conclusions are finally presented in Section~\ref{concl-s}.

{\bf Notation.} Throughout the paper, $\|v\|$ denotes the standard
Euclidean norm of a vector $v\in \Re^n$. For a symmetric tensor $S$ of order
$p$, $S[v_1, \ldots, v_p]$ is the result of applying $S$ to the vectors $v_1,
\ldots, v_p$, $S[v]^p$ is the result of applying $S$ to $p$ copies of the
vector $v$ and
\beqn{Tnorm}
\|S\|_{[p]} \eqdef \max_{\|v\|=1}  | S [v]^p |
= \max_{\|v_1\|= \cdots= \|v_p\|=1} | S[v_1, \ldots, v_p] |
\eeqn
(where the second equality results from Theorem~2.1 in \cite{ZhanLingQi12}) is the
associated induced norm for such tensors. If $S_1$ and $S_2$ are tensors,
$S_1\otimes S_2$ is their tensor product and $S_1^{k\otimes}$ is the product
of $S_1$ $k$ times with itself. For a real, sufficiently differentiable
univariate function $f$, $f^{(i)}$ denotes its $i$-th derivative and $f^{(0)}$ is
a synonym for $f$. For an integer $k$ and a real $\beta \in (0,1]$, we define
$(k+\beta)! \eqdef \prod_{\ell=1}^k(\beta+\ell)$ (this coincides with the
  standard factorial if $\beta=1$). As is usual, we also
define $0!=1$. If $M$ is a symmetric matrix, $\lambda_{\min}(M)$ is its
left-most eigenvalue. If $\alpha$ is a real, $\lceil \alpha \rceil$ and
$\lfloor \alpha \rfloor$ denote the smallest integer not smaller than $\alpha$
and the largest integer not exceeding $\alpha$, respectively. Finally
$\globmin_{x \in \calS}f(x)$ denotes the smallest value of $f(x)$ over $x \in
\calS$.

\numsection{High-order necessary conditions for optimality and the AR$p$ algorithm}
\label{arp-s}

Given $p \geq 1$, this paper considers the set-constrained optimization problem
\beqn{problem}
\min_{x \in \calF}  f(x),
\eeqn
where we assume that $\calF \subseteq \Re^n$ is closed and nonempty, and where
$f\in  \mathcal{C}^{p,\beta}(\Re^n)$, namely, that:
\begin{itemize}
\item  $f$ is
$p$-times continuously differentiable,

\item $f$ is bounded below by $f_{\rm low}$, and

\item the $p$-th derivative tensor of $f$ at $x$ is globally H\"{o}lder continuous, that is, 
there exist constants $L \geq 0$ and $\beta \in (0,1]$ such
that, for all $x,y \in \Re^n$,
\beqn{tensor-Hol}
\| \nabla_x^pf(x) - \nabla_x^pf(y) \|_{[p]} \leq L \| x-y \|^\beta.
\eeqn
\end{itemize}
Observe that convexity or even connectedness of $\calF$ is not requested.
Observe also that the more usual
case of \emph{Lipschitz continuous $p$-th derivative corresponds to
$\beta=1$}. We note that our assumption covers the continuous range of
objective function's smoothness from H\"{o}lder continuous gradients to Lipschitz
continuous $p$-th derivatives. In what follows, we assume that $\beta$ is known.

\noindent
If $T_p(x,s)$ is the standard $p$-th degree Taylor's expansion of $f$ about
$x$ computed for the increment $s$, that is
\beqn{Taylor-def}
T_p(x,s) \eqdef f(x) + \sum_{\ell=1}^p \frac{1}{\ell!}\nabla_x^\ell f(x)[s]^\ell,
\eeqn
\req{tensor-Hol}  provides crucial approximation bounds, whose proof can
be found in the appendix.

\llem{taylor-bounds-lemma}{
Let $f \in C^{p,\beta}(\Re^n)$, and $T_p(x,s)$ be the
Taylor approximation of $f(x+s)$ about $x$ given by \req{Taylor-def}.
Then for all $x,s \in \Re^n$,
\beqn{resf}
f(x+s) \leq T_p(x,s) + \bigfrac{L}{(p+\beta)!} \, \|s\|^{p+\beta},
\eeqn
\beqn{resder}
\| \nabla^j_x f(x+s) -  \nabla^j_s T_p(x,s) \|_{[j]}
\leq\bigfrac{L}{(p-j+\beta)!} \|s\|^{p-j+\beta}.
\ms (j = 1,\ldots, p).
\eeqn
}

In order to characterize minimizers of \req{problem}, we follow
\cite{CartGoulToin17c} and introduce, for given $\delta \in (0,1]$ and
$j\leq p$,
\beqn{phidef}
\phi_{f,j}^\delta(x)
\eqdef f(x)-\globmin_{\stackrel{x+d\in \calF}{\|d\|\leq\delta}}T_j(x,d),
\eeqn
which can be interpreted as the \emph{magnitude of the largest decrease achievable
on the Taylor's expansion of degree $j$ within the intersection a ball of radius
$\delta$ with the feasible set}. It was shown in \cite{CartGoulToin17c} that
$\phi_{f,j}^\delta(x)$ is a proper generalization of well-known unconstrained
optimality measures for low orders, in that, for $\delta = 1$,
\beqn{phi1-def}
\phi_{f,1}^\delta(x) = \|\nabla_x^1f(x)\|\,\delta,
\eeqn
\beqn{phi2-def}
\phi_{f,2}^\delta(x) = \left| \min[ 0, \lambda_{\min}\left(\nabla_x^2f(x)\right) \right|\,\delta^2
\eeqn
provided $\nabla_x^1f(x) = 0$, and also, if additionally $\nabla_x^2f(x)$ is positive
semi-definite, that
\beqn{phi3-def}
\phi_{f,3}^\delta = \|\;\mbox{projection of}\; \nabla_x^3f(x) \tim{onto the
  nullspace of} \nabla_x^2f(x) \;\|\,\delta^3.
\eeqn
At variance with other optimality measures, $\phi_{j,f}^\delta(x)$ is
well-defined for any order $j \geq 1$ and varies continuously when $x$ varies
continuoulsy in $\calF$. The role of the ``optimality radius'' $\delta$ in
\req{phidef} merits some discussion.
While the choice of $\delta=1$  is
adequate for retrieving known optimality conditions in the unconstrained case
for $j=1$, $j=2$ provided $\nabla_x^1f(x) = 0$, and $j=3$ provided additionnaly
$\nabla^2_xf(x)$ is positive semi-definite (as we have just seen), $\delta$
becomes important in other cases. Corollary~3.6 in \cite{CartGoulToin17c} indicates
that, when $\calF$ is convex, $q$-th-order necessary ``path-based'' optimality
conditions hold if
\beqn{limphi}
\lim_{\delta \rightarrow 0}  \frac{\phi_{f,j}^\delta(x)}{\delta^j} = 0
\tim{ for } j = 1,\ldots,q.
\eeqn
The limit for $\delta \rightarrow 0$ is necessary to capture the notion of
local minimizer for \req{problem}.  However, considering
$\phi_{f,j}^\delta(x)$ for non-vanishing $\delta$ has substantial advantages
from the point of view of optimization: while it may fail to indicate that $x$
is a local minimizer, it does so only by providing a direction leading to
values of $f$ below $f(x)$, thereby helping to avoid local but non-global
approximate solutions. We refer the reader to \cite{CartGoulToin17c} for a
further discussion, but conclude that considering fixed $\delta$ has strong
advantages when solving \req{problem}.

A special case is when $x$ is an isolated feasible point, that is a point
which is the sole intersection between $\calF$ and any sufficiently small
neighbourhood of $x$.  Such a point is clearly a local minimizer, and this is
reflected by the fact that $\phi_{f,q}^\delta(x) = 0$ for any $f$, any $q$ and
any sufficiently small $\delta$.

The main drawback of using $\phi_{f,j}^\delta(x)$ is, of course, that its
computation requires the global minimization of $T_p(x,d)$ in the intersection
of the ball of radius $\delta$ with $\calF$. We are not aware of an easy way to
do this in general\footnote{A small value of $\delta$ might help, but this
  computation remains NP-hard in most cases.} when
$n > 1$, which is why our analysis remains of an essentially theoretical
nature, as was the case for \cite{CartGoulToin17c}.  Note however that, albeit
potentially very difficult, solving this global minimization problem does not involve
calculating the value of $f$ or of any of its derivatives. In that sense, this
drawback is thus irrelevant for the worst-case evaluation complexity which
solely focuses on these evaluations.

Observe now that, if we were to relax the first-order condition
$\nabla_x^1f(x) = 0$ for unconstrained problems to $\|\nabla_x^1f(x)\|\leq
\epsilon$ and, at the same time, relax the second-order condition to
$\left| \min[ 0, \lambda_{\min}\left(\nabla_x^2f(x)\right) \right| \leq
  \epsilon$, we then deduce that
\beqn{term-q-2}
\phi_{f,2}^\delta(x) \leq \epsilon \delta + \half \epsilon \delta^2
= \epsilon \sum_{\ell=1}^2 \frac{\delta^\ell}{\ell!}.
\eeqn
A natural generalization of this observation is to define an
\emph{$(\epsilon,\delta)$-approximate $q$-th-order-necessary minimizer}
of $f$ as a point $x$ such that 
\beqn{term-q}
\phi_{f,q}^\delta(x) \leq \epsilon \chi_q(\delta)
\eeqn
where
\beqn{chidef}
\chi_q(\delta) \eqdef \sum_{\ell=1}^q \frac{\delta^\ell}{\ell!}.
\eeqn
Because \req{term-q} is a new way to look at approximate optimality and is
crucial for the rest of this paper, it is worthwhile to motivate and discuss it further.
\begin{enumerate}
\item When $\epsilon= 0$, \req{term-q} implies that the complicated
  path-based necessary optimality conditions derived in \cite{CartGoulToin17c} do hold.
  This results from the fact that these latter conditions merely express that
  the Taylor's model of order $q$ cannot decrease close enough to $x$ along
  any feasible polynomial path emanating from $x$, which is clearly the case
  if $x$ is a global minimizer of the same models in the intersection of the
  feasible set and a ball of radius $\delta$ centered at $x$. By continuity, these path-based
  conditions must therefore hold in the limit under \req{term-q} when $\epsilon$ tends
  to zero.  The role of \req{term-q} as a condition for approximate
  minimization is thus coherent and consistent with known necessary conditions.
\item Inspired by \req{limphi}, the stronger approximate optimality
  condition
  \beqn{too-strong}
  \phi_{f,j}^\delta(x) \leq \epsilon \, \delta^j \tim{ for }j\in \ii{q}
  \eeqn
  was used in \cite{CartGoulToin17c} instead of \req{term-q}. Our main reason
  to prefer \req{term-q} is the following. Observe that \req{too-strong}
  implies in particular that $\phi_{f,q}^\delta(x) \leq \epsilon \delta^q$,
  which in turn implies, for $\delta$ small enough for the first-order term to
  dominate, that $\phi_{f,1}^\delta(x) \leq \epsilon \delta^q$.  In the
  unconstrained case (for example), this requires $\|\nabla_x^1f(x_k)\|\leq
  \epsilon \delta^{q-1}$, imposing an inordinate level of first-order
  optimality, much stronger than the standard condition
  $\|\nabla_x^1f(x_k)\|\leq \epsilon$. No such difficulty arises with
  \req{term-q} because the right-hand side of the condition involves all
  powers of $\delta$, which is not the case of the right-hand side of
  \req{too-strong}. Note however that the vital continuity properties of
  $\phi_{f,q}^\delta$ are not affected by the choice of the right-hand side,
  and are thus inherited by \req{term-q}.
\item For given $\delta \in (0,1]$, \req{term-q} does not imply that
  $\phi_{f,j}^\delta(x) \leq \epsilon \chi_j(\delta)$ for $j \in \ii{q-1}$,
  although the violation of this condition tends to zero with
  $\delta$\footnote{\label{gap}When $\delta$ tends to zero,
    the terms of orders $j+1$ and higher in the Taylor's expansion defining $\phi_{f,q}^\delta(x)$
    and $\chi_q(\delta)$ become negligible
    compared to the first $j$.}. This slight blemish can be cured by requiring
  that $\phi_{f,j}^\delta(x) \leq \epsilon \chi_j(\delta)$ for $j\in \ii{q}$
  instead of \req{term-q}, but we claim that the benefit of this stronger
  definition is outweighted by the need to perform $q-1$ additional constrained
  global minimizations, and therefore focus our exposition to the case using the
  simpler \req{term-q}.
\end{enumerate}

\noindent
In order to further justify \req{term-q}, we now make more explicit the
``minimizing guarantees'' provided by this approximate optimality condition, by
formulating a result analogous to Theorem~3.7 in \cite{CartGoulToin17c}.  This
result gives a lower bound on the value of $f(x)$ in the feasible
neighbourhood of an $(\epsilon,\delta)$-approximate $q$-th-order-necessary minimizer.

\lthm{analog-3.7}{Suppose that $f$ is $p$ times continuouly differentiable
and that $\nabla^q_x f$ is $\beta$-H\"{o}lder continuous with constant $L$ (in
the sense of \req{tensor-Hol} with $p=q$) in an open neighbourhood of radius
$\delta \in (0,1]$ of some $x \in \calF$. Suppose also that $x$ is an
$(\epsilon,\delta)$-approximate $q$-th-order-necessary minimizer of $f$ in the
sense of \req{term-q}.  Then 
\beqn{near-minimizer}
f(x+d) \geq f(x) - 2 \epsilon \chi_q(\delta)
\tim{ for all $d$ with }
x+d \in \calF
\tim{and}
\|d\| \leq \min\left[ \delta, \left(\frac{(q+1)!\epsilon}{L}\right)^{\frac{1}{q+\beta-1}}\right].
\eeqn
}

\proof{Using the triangle inequality, \req{tensor-Hol}, \req{resf} and
\req{term-q}, we obtain that
\[
\begin{array}{lcl}
f(x+d)
& \geq & f(x+d)-T_q(x,d) + T_q(x,d) \\*[1.5ex]
& \geq & -|f(x+d)-T_q(x,d)|+ T_q(x,0)-\phi_{f,q}^\delta(x) \\*[1.5ex]
& \geq & -\bigfrac{L}{(q+1)!}\|d\|^{q+\beta}+f(x) - \epsilon \chi_q(\delta).
\end{array}
\]
Thus, if $\|d\|\leq \delta$,
\[
f(x+d) \geq f(x) - \bigfrac{L}{(q+1)!}\|d\|^{q+\beta-1}\,\delta - \epsilon
\chi_q(\delta)
\]
and the desired bound follows from the fact that $\delta \leq \chi_q(\delta)$.
} 

\noindent
In order to find $(\epsilon,\delta)$-approximate $q$-th-order-necessary
minimizers, we consider applying a variant of the AR$p$ algorithm to
\req{problem}. This algorithm, described as Algorithm~\ref{algo}
\vpageref{algo}, is of the regularization type in that, at each iterate $x_k$,
a step $s_k$ is computed which approximately minimizes (in a sense defined
below) the model
\beqn{model}
m_k(s) = T_p(x_k,s) + \frac{\sigma_k}{(p+\beta)!} \|s\|^{p+\beta}
\eeqn
subject to $x_k+s\in \calF$,
where $p$ in an integer such that $p\geq q$ and $\sigma_k \geq \sigma_{\min}$
is a ``regularization parameter''. 

\algo{algo}{AR$p$ for $(\epsilon,\delta)$-approximate
  $q$-th-order-necessary minimizers}
{
\vspace*{-0.3 cm}
\begin{description}
\item[Step 0: Initialization.]
  An initial point $x_0\in\calF$  and an initial regularization parameter $\sigma_0>0$
  are given, as well as an accuracy level  $\epsilon \in (0,1)$.  The
  constants $\delta_{-1}$, $\varpi$, $\theta$, $\eta_1$, $\eta_2$, $\gamma_1$, $\gamma_2$,
  $\gamma_3$ and $\sigma_{\min}$ are also given and satisfy
  \beqn{eta-gamma2}
  \begin{array}{c}
\varpi \in (0,1], \;\; \theta > 0,  \;\; \delta_{-1} \in (0,1], \;\;
    \sigma_{\min} \in (0, \sigma_0], \;\;
      0 < \eta_1 \leq \eta_2 < 1 \\
      \tim{and} 0< \gamma_1 < 1 < \gamma_2 < \gamma_3.
\end{array}
\eeqn
Compute $f(x_0)$ and set $k=0$.

\item[Step 1: Test for termination. ]
  Evaluate $\{\nabla^i_x f(x_k)\}_{i=1}^q$.
  If \req{term-q} holds with $\delta = \delta_{k-1}$, terminate with the approximate
  solution $x_\epsilon=x_k$. Otherwise compute $\{\nabla^i_x f(x_k)\}_{i=q+1}^p$.

\item[Step 2: Step calculation. ] Attempt to compute a step $s_k$ such
that $x_k+s_k\in \calF$ and an optimality radius $\delta_k\in (0,1]$ by
approximately minimizing the model $m_k(s)$ in
the sense that
\vspace*{-1mm}
\beqn{descent2}
m_k(s_k) < m_k(0)
\vspace*{-1mm}
\eeqn
and either
\vspace*{-1mm}
\beqn{mterm1}
\|s_k\| \geq \varpi \epsilon^{\frac{1}{p-q+\beta}}
\eeqn
or
\beqn{mterm2}
\phi_{m_k,q}^{\delta_k}(s_k) \leq \frac{\theta\|s_k\|^{p-q+\beta}}{(p-q+\beta)!}\,\chi_q(\delta_k).
\eeqn
If no such step exist, terminate with the approximate solution
$x_\epsilon=x_k$.

\item[Step 3: Acceptance of the trial point. ]
Compute $f(x_k+s_k)$ and define
\beqn{rhokdef2}
\rho_k = \frac{f(x_k) - f(x_k+ s_k)}{T_p(x_k,0) - T_p(x_k,s_k)}.
\eeqn
If $\rho_k \geq \eta_1$, then define
$x_{k+1} = x_k + s_k$; otherwise define $x_{k+1} = x_k$.

\item[Step 4: Regularization parameter update. ]
Set
\beqn{sigupdate2}
\sigma_{k+1} \in \left\{ \begin{array}{ll}
{}[\max(\sigma_{\min}, \gamma_1\sigma_k), \sigma_k ]  & \tim{if} \rho_k \geq \eta_2, \\
{}[\sigma_k, \gamma_2 \sigma_k ]          &\tim{if} \rho_k \in [\eta_1,\eta_2),\\
{}[\gamma_2 \sigma_k, \gamma_3 \sigma_k ] & \tim{if} \rho_k < \eta_1.
  \end{array} \right.
\eeqn
Increment $k$ by one and go to Step~1 if $\rho_k\geq \eta_1$, or to Step~2 otherwise.
\end{description}
}

\noindent
A few comments are useful at this stage.
\begin{enumerate}
\item Since $\sigma_k \geq \sigma_{\min}$ by \req{sigupdate2}, we have that
  $m_k(s)$ is bounded below as a function of $s$ and the existence of a
  constrained global minimizer $s_k^*$ is guaranteed.
\item Step~2 requires, that, for $s_k\neq 0$, we also compute $\delta_k$.  This
  is easy for orders one and two.  If $q=1$, the formula for a global
  minimizer $s_k^*$ is analytic and $\delta_k=1$ is always acceptable.  The
  situation is similar for $q= 2$, where $s_k^*$ can be assessed using a
  trust-region method whose radius is $\delta_k=1$ (more details are provided
  at the end of Section~\ref{upperbound-s}). The task is more difficult for
  higher orders where one may have to rely on the arguments of
  Lemma~\ref{step-ok-l} below, or use different subproblems with decreasing
  values of $\delta$.  However, none of these computations involve the
  evaluation of $f$ or its derivatives, and therefore the evaluation
  complexity bound discussed in this paper is unaffected.
\item That one needs to consider the second case in Step~2 (where no step exists
  satisfying \req{descent2} -- \req{mterm2}) can be seen by examining the
  following one-dimensional example. Let $p=q=3$ and $\beta=1$,  and
  suppose that $\delta_{k-1}=1$, $T_q(x_k,s) = s^2-2s^3$ and $\sigma_k = 4!=24$.
  Then $m_k(s) = s^2-2s^3 + s^4= s^2(s-1)^2$ and the origin is a
  global minimizer of the model (and a local minimizer of $T_q(x_k,s)$) but yet
  $T_q(x_k,\delta) = -1$, yielding that $\phi_{f,q}^{\delta_{k-1}}(x_k) = 1 >
  \epsilon \chi_q(1)$ for $\epsilon\leq 1/\chi_q(1)=\sfrac{4}{7}$. Thus,
  Step~1 with $\delta_{k-1}=1$ has failed to identify that termination was
  possible. In addition, we see that, at variance with the cases $q=1$ and
  $q=2$, a global minimizer of the model \req{model} may not, for $q\geq 3$,
  be a global minimizer of its $q$-th order Taylor's expansion in the
  intersection of $\calF$ and a ball of arbitrary radius: we may have 
  to restrict this radius (to $\delta_{k-1}=\half$ in our example) for this
  important property to hold (see Lemma~\ref{step-ok-l} below).
\item If \req{mterm1} holds, the  possibly expensive computation of
  $\phi_{m_k,q}^{\delta_k}(s_k)$ in \req{mterm2} is
  unnecessary and $\delta_k$ may be chosen arbitrarily in $(0,1]$.
\item We assume the availability of a feasible starting point, which is
  without loss of generality for inexpensive constraints.
\item Before termination, each successful iteration requires the evaluation
  of $f$ and its first $p$ derivative tensors, while only the evaluation of
  $f$ is needed at unsuccessful ones.
\item The mechanism of the algorithm ensures the non-increasing nature of the
  sequence $\{f(x_k)\}_{k\geq 0}$.
\end{enumerate}
Iterations for which $\rho_k \geq \eta_1$
(and hence $x_{k+1}=x_k+s_k$) are called ``successful'' and we denote by $\calS_k
\eqdef \{ 0 \leq j \leq k \mid \rho_j \geq \eta_1 \}$ the index set of all
successful iterations between 0 and $k$.
We immediately observe that the total number of iterations (successful or not)
can be bounded as a function of the number of successful ones (and include a
proof in the appendix).

\llem{SvsU}{
\cite[Theorem~2.4]{BirgGardMartSantToin17}
The mechanism of Algorithm~\ref{algo} guarantees that, if
\beqn{sigmax}
\sigma_{k} \leq \sigma_{\max},
\eeqn
for some $\sigma_{\max} > 0$, then
\beqn{unsucc-neg}
k +1 \leq |\calS_k| \left(1+\frac{|\log\gamma_1|}{\log\gamma_2}\right)+
\frac{1}{\log\gamma_2}\log\left(\frac{\sigma_{\max}}{\sigma_0}\right).
\eeqn
}

We also verify that the algorithm is well-defined in the sense that either a step
$s_k$ satisfying \req{descent2}--\req{mterm2} can always be found, or
termination is justified.  For
unconstrained problems with $q\in\{1,2\}$, the first possibility directly results from the
observation that $\phi_{m_k,j}^\delta(s_k)$ (as given by \req{phi1-def}-\req{phi3-def}
for $f=m_k$ and $j\in\{1,2,3\}$) can be made suitably small at a global minimizer of the
model. The situation is more complicated for other cases.  In order to clarify
it, we first state a useful technical lemma, whose proof is in the appendix.

\llem{njsp-l-a}{
Let $s$ be a vector of $\Re^n$.  Then
\beqn{npsp-a}
\|\, \nabla_s^j \big(\|s\|^{p+\beta} \big) \, \|_{[j]}
= \frac{(p+\beta)!}{(p-j+\beta)!}\|s\|^{p-j+\beta}
\tim{for} j\in \iibe{0}{p}
\eeqn
and
\beqn{npsp-b}
\|\, \nabla_s^{p+1} \big(\|s\|^{p+\beta} \big) \, \|_{[p+1]}
= \beta \,(p+\beta)! \,\|s\|^{\beta-1}.
\eeqn
}

We now provide reasonable sufficient conditions for a nonzero step $s_k$ and
an optimality radius $\delta_k$ to satisfy \req{descent2}--\req{mterm2}.

\llem{step-ok-l}{Suppose that $s_k^*$ is a global minimizer of $m_k(s)$
  under the constraint that $x_k+s \in \calF$, such $m_k(s_k^*)< m_k(0)$.
  Then there exist a neighbourhood of $s_k^*$ and a range of sufficiently
  small $\delta$ such that \req{descent2} and \req{mterm2}
  hold for any $s_k$ in the intersection of this neighbourhood with $\calF$
  and any $\delta_k$ in this range.}

\proof{
Let $s_k^*$ be the global minimizer of the model $m_k(s)$ over all $s$ such
that $x_k+s \in \calF$. Since $m_k(s_k^*)< m_k(0)$, we have that $s_k^* \neq
0$. By Taylor's theorem, we have that, for all $d$,
\[
0 \leq  m_k(s_k^*+d)-m_k(s_k^*)
= \sum_{\ell=1}^p \frac{1}{\ell!}\nabla_s^\ell m_k(s_k^*)[d]^\ell
+ \frac{1}{(p+1)!} \nabla_s^{p+1}m_k(s_k^*+\xi d)[d]^{p+1}
\]
for some $\xi \in (0,1)$. Thus, using the triangle inequality, \req{model} and \req{npsp-b},
\beqn{sl-1}
\begin{array}{lcl}
- \bigsum_{\ell=1}^q \bigfrac{1}{\ell!}\nabla_s^\ell m_k(s_k^*)[d]^\ell
& \leq & \bigsum_{\ell=q+1}^p \bigfrac{\|d\|^\ell}{\ell!}\|\nabla_s^\ell m_k(s_k^*)\|_{[\ell]}
   + \bigfrac{\|d\|^{p+1}}{(p+1)!} \|\nabla_s^{p+1}m_k(s_k^*+\xi d)\|_{[p+1]}\\*[2ex]
& = & \bigsum_{\ell=q+1}^p \bigfrac{\|d\|^\ell}{\ell!}\|\nabla_s^\ell m_k(s_k^*)\|_{[\ell]}
   + \beta \sigma_k \bigfrac{\|d\|^{p+1}}{(p+1)!}\|s_k^*+\xi d \|^{\beta-1}.
\end{array}
\eeqn
Since $s_k^* \neq 0$, we may then choose $\delta_k < \|s_k^*\|$ such that, for
every $d$ with $\|d\| \leq \delta_k$, $\|s_k^*+\xi d \|\geq \half \|s_k^*\|>0$ and
\beqn{sl-2}
\bigsum_{\ell=q+1}^p \bigfrac{\|d\|^\ell}{\ell!}\|\nabla_s^\ell m_k(s_k^*)\|_{[\ell]}
+ 2^{1-\beta}\beta \sigma_k\bigfrac{\|d\|^{p+1}}{(p+1)!}\|s_k^* \|^{\beta-1}
\leq \frac{\theta \|s_k^*\|^{p-q+\beta}}{2 (p-q+\beta)!}\,\|d\|.
\eeqn
Hence we deduce from \req{sl-1} and \req{sl-2} that, for $\|d\|\leq\delta_k$,
\[
-\sum_{\ell=1}^q\frac{1}{\ell!}\nabla_s^\ell m_k(s_k^*)[d]^\ell
\leq \frac{\theta \|s_k^*\|^{p-q+\beta}}{2(p-q+\beta)!}\,\delta_k
\leq \frac{\theta \|s_k^*\|^{p-q+\beta}}{2(p-q+\beta)!}\,\chi_q(\delta_k),
\]
where the last inequality follows from \req{chidef}.
Continuity of $m_k$ and its derivatives and the inequality $m_k(s_k^*) <
m_k(0)$ then imply that there exists a neighbourhood of $s_k^*\neq 0$ such
that \req{descent2} holds and 
\[
-\sum_{\ell=1}^q\frac{1}{\ell!}\nabla_s^\ell m_k(s)[d]^\ell
\leq \frac{\theta \|s\|^{p-q+\beta}}{(p-q+\beta)!}\,\chi_q(\delta_k).
\]
for all $s$ in this neighbourhood and all $d$ with $\|d\|\leq \delta_k$.
This yields that, for all such $s$ with $x_k+s \in \calF$,
\[
\phi_{m_k,q}^{\delta_k}(s)
= \max\Bigg[
      0,
      \globmax_{\mystack{\|d\|\leq\delta_k}{x_k+d\in\calF}}
      \left(-\sum_{\ell=1}^q\frac{1}{\ell!}\nabla_s^\ell m_k(s_k)[d]^\ell\right)
      \Bigg]
\leq \frac{\theta \|s\|^{p-q+\beta}}{(p-q+\beta)!}\,\chi_q(\delta_k),
\]
as requested.
} 

\noindent
As can be seen in the proof of this lemma, $\delta_k$ may need to be small if
any of the tensors
\[
\nabla_s^\ell m_k(s_k^*)
= \sum_{j=\ell}^p \frac{1}{j!}\nabla_s^j m_k(0)[s_k^*]^{j-\ell}
\]
for $\ell \in \ii{p+1}$ has a large norm. This may occur in particular if
$\beta$ and $\|s_k^*\|$ are both close to zero, as is shown by the last term
in the left-hand side of \req{sl-2}.  We also note that 
\req{mterm2} obviously holds for $s_k=s_k^*$ if $x_k+s_k^*$ is an isolated
feasible point. It now remains to verify that it is justified to terminate in 
Step~2 when no suitable nonzero step can be found.

\llem{term-step2-l}{Suppose that the algorithm terminates in Step~2 of
  iteration $k$ with $x_\epsilon=x_k$. Then there exists a $\delta \in
  (0,1]$ such that \req{term-q} holds for $x = x_\epsilon$ and $x_\epsilon$ is
  an $(\epsilon,\delta)$-approximate $q$th-order-necessary minimizer.
}

\proof{ Given Lemma~\ref{step-ok-l}, if the algorithm terminates within
  Step~2, it must be because every global minimizer $s_k^*$ of $m_k(s)$ under
  the constraints $x_k+s\in \calF$ is such that $m_k(s_k^*)\geq m_k(0)$. In
  that case, $s_k^*=0$ is one such global minimizer and we have that, for all
  $d$,
  \[
  0 \leq  m_k(d) - m_k(0)
  = \sum_{\ell=1}^q\frac{1}{\ell!}\nabla_x^jf(x_k)[d]^j
  + \sum_{\ell=q+1}^p\frac{1}{\ell!}\nabla_x^jf(x_k)[d]^j
  + \frac{\sigma_k}{(p+\beta)!}\|d\|^{p+\beta}.
  \]
  We may now choose $\delta \in (0,1]$ small enough to ensure that, for all
  $d$ with $\|d\|\leq \delta$,
  \beqn{sdelta}
  \left|\sum_{\ell=q+1}^p\frac{1}{\ell!}\nabla_x^jf(x_k)[d]^j
  + \frac{\sigma_k}{(p+\beta)!}\|d\|^{p+\beta}\right|
  \leq \epsilon \|d\|
  \leq \epsilon \, \chi_q(\delta),
  \eeqn
  which in turn implies that, for all $d$ with $\|d\|\leq \delta$,
  \[
  \phi_{f,q}^{\delta}(x_k)
  = \max\Bigg[
      0,
      \globmax_{\mystack{\|d\|\leq\delta}{x_k+d\in\calF}}
      \left(-\sum_{\ell=1}^q\frac{1}{\ell!}\nabla_x^\ell f(x_k)[d]^\ell\right)
      \Bigg]
  \leq \epsilon\,\chi_q(\delta_k),
  \]
  concluding the proof.
} 

\noindent
Observe that, in this proof, we could have chosen $\delta$ small enough to
ensure
\[
\frac{\sigma_k}{(p+\beta)!}\|d\|^{p+\beta} \leq \epsilon \chi_p(\delta)
\]
instead of \req{sdelta}, yielding 
$\phi_{f,p}^\delta(x_k) \leq \epsilon \chi_p(\delta)$, which is a stronger
necessary optimality condition than \req{term-q}.
Together, Lemmas~\ref{step-ok-l} and \ref{term-step2-l} ensure that
Algorithm~\ref{algo} is well-defined.

\numsection{An upper bound on the evaluation complexity}
\label{upperbound-s}

The proofs of the following two lemmas are very similar to corresponding
results in \cite{BirgGardMartSantToin17} and hence we again defer them to the
appendix (but still include them for completeness, as the algorithm has
changed).

\llem{Dm-lemma}{
The mechanism of Algorithm~\ref{algo} guarantees that, for all $k  \geq 0$,
\beqn{Dphi}
T_p(x_k,0) - T_p(x_k,s_k) \geq \frac{\sigma_k}{(p+\beta)!} \|s_k\|^{p+\beta},
\eeqn
and so \req{rhokdef2} is well-defined.
}

\llem{sigmaupper-lemma}{
Let $f\in \mathcal{C}^{p,\beta}(\Re^n)$.
Then, for all $k\geq 0$,
\beqn{sigmaupper}
\sigma_k
\leq \sigma_{\max} \eqdef \max\left[ \sigma_0,\frac{\gamma_3 L}{1-\eta_2}\right].
\eeqn
}

\noindent
We are now in position to prove the crucial lower bound on the
step length.

\llem{longs-j-lemma}{
Let $f\in \mathcal{C}^{p,\beta}(\Re^n)$. Then, for all  $k\geq 0$ such that
Algorithm~\ref{algo} does not terminate at iteration $k+1$,
\beqn{longs-j}
\|s_k\|\geq \kappa_s \epsilon^{\frac{1}{p-q+\beta}},
\eeqn
where
\beqn{kappas-def}
\kappa_s \eqdef
\min\left[\varpi,\left(\frac{(p-q+\beta)!}
 {(L+\sigma_{\max}+\theta)}\right)^{\frac{1}{p-q+\beta}}\right].
\eeqn
}

\proof{
If $\|s_k\| > \varpi \epsilon^{\frac{1}{p-q+\beta}}$, the result is
obvious.  Suppose now that $\|s_k\| \leq \varpi \epsilon^{\frac{1}{p-q+\beta}}$.
Since the algorithm does not terminate at iteration $k+1$, we have that
\beqn{big-phi-plus}
\phi_{f,q}^{\delta_k}(x_{k+1}) > \epsilon \chi_q(\delta_k)
\eeqn
Let the global minimum in the definition of $\phi_{f,q}^{\delta_k}(x_{k+1})$ be
achieved at $d$ with $\|d\|\leq \delta_k$. Since
$\phi_{f,q}^{\delta_k}(x_{k+1}) >0$, we have from \req{phidef} that 
\[
\bigsum_{\ell=1}^q \bigfrac{1}{\ell!}\nabla^\ell_xf(x_{k+1})[d]^\ell < 0
\]
Then, successively using \req{phidef} for
$f$ at $x_{k+1}$, the triangle inequality, \req{model}, \req{Tnorm}
and \req{npsp-a}, we deduce that
\beqn{big-bound1}
\begin{array}{lcl}  
\phi_{f,q}^{\delta_k}(x_{k+1})
&   =  & - \bigsum_{\ell=1}^q \bigfrac{1}{\ell!}\nabla^\ell_xf(x_{k+1})[d]^\ell \\*[2ex]
&   =  & - \bigsum_{\ell=1}^q \bigfrac{1}{\ell!}\nabla^\ell_xf(x_{k+1})[d]^\ell
         + \bigsum_{\ell=1}^q \bigfrac{1}{\ell!}\nabla^\ell_sT_p(x_k,s_k)[d]^\ell
         - \bigsum_{\ell=1}^q \bigfrac{1}{\ell!}\nabla^\ell_sT_p(x_k,s_k)[d]^\ell\\*[2ex]
&      & - \bigfrac{\sigma_k}{(p+\beta)!}
           \bigsum_{k=\ell}^q\bigfrac{1}{\ell!}\left(\nabla^\ell_s \big[\|s\|^{p+\beta}\big](s_k)\right)[d]^\ell
         + \bigfrac{\sigma_k}{(p+\beta)!}
           \bigsum_{k=\ell}^q\bigfrac{1}{\ell!}\left(\nabla^\ell_s \big[\|s\|^{p+\beta}\big](s_k)\right)[d]^\ell\\*[2ex]
& \leq &  \left| \bigsum_{\ell=1}^q \bigfrac{1}{\ell!}\Big[\nabla^\ell_xf(x_{k+1})-\nabla^\ell_sT_p(x_k,s_k)\Big][d]^\ell\right| \\*[2ex]
&      & \hspace*{2cm} - \bigsum_{\ell=1}^q \bigfrac{1}{\ell!}
         \left(\nabla^\ell_s\left[ T_p(x_k,s)+\bigfrac{\sigma_k}{(p+\beta)!}\|s\|^{p+\beta}\right]_{s=s_k}\right)[d]^\ell \\*[2ex]
&      & \hspace*{2cm} + \bigfrac{\sigma_k}{(p+\beta)!}\left|
         \bigsum_{k=\ell}^q\bigfrac{1}{\ell!}\left(\nabla^\ell_s \big[\|s\|^{p+\beta}\big]_{s=s_k}\right)[d]^\ell\right|\\*[2ex]
& \leq &  \bigsum_{\ell=1}^q\bigfrac{L}{\ell!(p-\ell+\beta)!}\|s_k\|^{p-\ell+\beta}\delta_k^\ell \\*[2ex]
&      & \hspace*{2cm} - \bigsum_{\ell=1}^q \bigfrac{1}{\ell!} \nabla^\ell_s m_k(s_k)[d]^\ell 
          + \bigsum_{\ell=1}^q \bigfrac{\sigma_k}{\ell!(p-\ell+\beta)!}\|s_k\|^{p-\ell+\beta}\delta_k^\ell\\*[2ex]
\end{array}         
\eeqn
Now, since $\|d\| \leq \delta_k$, and using \req{phidef} for $m_k$ at $s_k$,
\[
- \bigsum_{\ell=1}^q \bigfrac{1}{\ell!} \nabla^\ell_s m_k(s_k)[d]^\ell
\leq \max\left[0,-\bigsum_{\ell=1}^q \bigfrac{1}{\ell!} \nabla^\ell_s
  m_k(s_k)[d]^\ell\right]
\leq \phi_{m_k,q}^{\delta_k}(s_k).
\]
Therefore, using \req{mterm2} and \req{big-bound1}, we have that
\beqn{big-bound2}
\begin{array}{lcl}  
\phi_{f,q}^{\delta_k}(x_{k+1})
& \leq & \bigsum_{\ell=1}^q \bigfrac{L}{\ell!(p-\ell+\beta)!}\|s_k\|^{p-\ell+\beta}\delta_k^\ell
         + \bigfrac{\theta\,\chi_q(\delta_k)}{(p-q+\beta)!} \|s_k\|^{p-q+\beta}\\*[2ex]
&      & \hspace*{2cm} + \bigsum_{\ell=1}^q \bigfrac{\sigma_k}{\ell!(p-\ell+\beta)!}\|s_k\|^{p-\ell+\beta}\delta_k^\ell \\*[2ex]
& \leq & \bigfrac{\Big[ L+\sigma_k + \theta\Big] \chi_q(\delta_k)}{(p-q+\beta)!}\|s_k\|^{p-q+\beta},
\end{array}         
\eeqn
where we have used the fact that $\|s_k\| \leq \varpi
\epsilon^{\frac{1}{p-q+\beta}} \leq 1$ to deduce the last inequality.
As a consequence, \req{big-phi-plus} implies that
\[
\|s_k\|
\geq
\left[\frac{\epsilon(p-q+\beta)!}{(L+\sigma_k+\theta)}\right]^{\frac{1}{p-q+\beta}}
\]
and \req{longs-j} then immediately follows from \req{sigmaupper}.
} 

\noindent
The bound given by this lemma is another indication that choosing $\theta$ of
the order of $L$ (when this is known a priori) makes sense.

We now combine all the above results to deduce an upper bound on the maximum
number of successful iterations, from which a final complexity bound
immediately follows.

\lthm{upper-theorem}{
Let $f\in \mathcal{C}^{p,\beta}(\Re^n)$.
Then, given $\epsilon \in (0,1)$,
Algorithm~\ref{algo} needs at most
\[
\left \lfloor \kappa_p ( f(x_0)- \flow)
\left( \epsilon^{-\frac{p+\beta}{p-q+\beta}} \right)
\right \rfloor+1
\]
successful iterations (each involving one evaluation of $f$ and its $p$ first derivatives)
and at most
\beqn{final-upper}
\left\lfloor \left \lfloor
\kappa_p ( f(x_0)- \flow)
\left( \epsilon^{-\frac{p+\beta}{p-q+\beta}} \right) +1
\right \rfloor
                 \left(1+\frac{|\log\gamma_1|}{\log\gamma_2}\right)+
\frac{1}{\log\gamma_2}\log\left(\frac{\sigma_{\max}}{\sigma_0}\right)\right\rfloor
\eeqn
iterations in total to produce an iterate $x_\epsilon$ such that
\req{term-q} holds, where $\sigma_{\max}$ is given by
\req{sigmaupper} and where 
\[
\kappa_p \eqdef \frac{(p+\beta)!}{\eta_1 \sigma_{\min}}
\max\left\{ \varpi^{-(p+\beta)},\left[\frac{(L+\sigma_{\max}+\theta)}{(p-q+\beta)!}\right]^{\frac{p+\beta}{p-q+\beta}}\right\}.
\]
}

\proof{At each successful iteration $k$ before termination, we have the guaranteed decrease
\beqn{fdec}
f(x_k)-f(x_{k+1})
\geq \eta_1 (T_p(x_k,0)-T_p(x_k,s_k))
\geq \bigfrac{\eta_1 \sigma_{\min}}{(p+\beta)!} \;\|s_k\|^{p+\beta}
\eeqn
where we used \req{rhokdef2}, \req{Dphi} and \req{sigupdate2}. Moreover  we
deduce from \req{fdec}, \req{longs-j} and \req{sigmaupper} that 
\beqn{eps1-decr}
f(x_k)-f(x_{k+1}) \geq \kappa_p^{-1} \epsilon_j^{\frac{p+\beta}{p-q+\beta}}
\tim{where}
\kappa_p^{-1} \eqdef
\bigfrac{\eta_1 \sigma_{\min}\kappa_s}{(p+\beta)!}.
\eeqn
Thus, since $\{f(x_k)\}$ decreases monotonically,
\[
f(x_0)-f(x_{k+1})
\geq \kappa_p^{-1} \,\epsilon_j^{\frac{p+\beta}{p-q+\beta}}\,|\calS_k|.
\]
Using that $f$ is bounded below by $f_{\rm low}$, we conclude 
\beqn{Sk1}
| \calS_k |
\leq \frac{f(x_0) - \flow}{\kappa_p^{-1}} \epsilon_j^{-\frac{p+\beta}{p-q+\beta}}
\eeqn
until termination. The desired bound on the number of
successful iterations follows from combining \req{Sk1}.
Lemma~\ref{SvsU} is then invoked to compute the upper bound on the total
number of iterations.
}

\noindent
In particular, if the $p$-th derivative of $f$ is assumed to be globally
Lipschitz rather than merely H\"{o}lder continuous (i.e.\ if $\beta = 1$), the
bound \req{final-upper} on the maximum number of evaluations becomes
\beqn{final-bound-Lip}
\left \lfloor
\left \lfloor
\kappa_p ( f(x_0)- \flow)
\left( \epsilon^{-\frac{p+1}{p-q+1}} \right)+1\right \rfloor
                 \left(1+\frac{|\log\gamma_1|}{\log\gamma_2}\right)+
\frac{1}{\log\gamma_2}\log\left(\frac{\sigma_{\max}}{\sigma_0}\right)\right \rfloor
\eeqn
where
\[
\kappa_p \eqdef \frac{(p+1)!}{\eta_1 \sigma_{\min}}
\max\left\{ \varpi^{p+\beta},\left[\frac{q!(L+\sigma_{\max}+\theta)(e-1)}{(p-q+1)!}\right]^{\frac{p+1}{p-q+1}}\right\}.
\]
This worst-case evaluation bound generalizes known bounds
for $q= 1$ (see \cite{BirgGardMartSantToin17}) or $q=2$ (see\cite{CartGoulToin17e}) and
significantly improve upon the bounds in $O(\epsilon^{-(q+1)})$ given by
\cite{CartGoulToin17c} for a more stringent termination rule. It also extends
the results obtained in \cite{CartGoulToin12b} for convexly-constrained
problems with $q=1$ by allowing the significantly broader class of inexpensive constraints.

We also note that it is possible to weaken the assumption that $\nabla_x^pf$
must satisfy the H\"{o}lder inequality \req{tensor-Hol} for every $x,y \in
\Re^n$ (as required in the beginning of Section~\ref{arp-s}). The weakest
possible smoothness assumption is to require that \req{tensor-Hol} holds only
for points belonging to the same segment of the ``path of iterates''
$\cup_{k\geq 0} [x_k,x_{k+1}]$ (this is necessary for the proof of
Lemma~\ref{taylor-bounds-lemma}). As this path joining feasible iterates may
be hard to predict a priori, one may instead use the monotonic character of
Algorithm~\ref{algo} and require \req{tensor-Hol} to hold for all $x,y$ in the
intersection of $\calF$ with the level set $\{x \in \Re^n \mid f(x)\leq
f(x_0)\}$. Again, it may be hard to determine this set and to ensure that it
contains the path of iterates, and one may then resort to requiring
\req{tensor-Hol} to hold in the whole of $\calF$, which must then be convex to
ensure the desired H\"{o}lder property on every segment $[x_k,x_{k+1}]$.

\numsection{Seeking $\epsilon$-approximate second-order-necessary minimizers}
\label{2ndorder-s}

We now discuss the particular and much-studied case
where second-order minimizers are sought for unconstrained problems with
Lipschitz continuous Hessians (that is $p\geq q=2$, $\calF=\Re^n$ and
$\beta=1$). As we now show, a specialization of Algorithm~\ref{algo} to this 
case is very close (but not identical) to well-known methods.
Let us consider Step~1 first.  The computation of
$\phi_{f,2}^{\delta_{k-1}}(x_k)$ then reduce to 
\beqn{phi-2}
\phi_{f,2}^{\delta_{k-1}}(x_k)
= \max\Bigg[
      0,
      - \globmin_{\|d\|\leq \delta_{k-1}}
      \Big(\nabla_x^1f(x_k)^Td + \half d^T\nabla_x^2f(x_k)d\Big)
      \Bigg],
\eeqn
which amounts to solving a standard trust-region subproblem with radius $\delta_{k-1}$ (see
\cite{ConnGoulToin00}). Hence verifying \req{phi-2} or testing the more usual
approximate second-order criterion
\beqn{strong-2}
\|\nabla_x^1f(x_k)\| \leq \epsilon
\tim{and}
\lambda_{\min}\Big(\nabla_x^2f(x_k)\Big) \geq -\epsilon,
\eeqn
have very similar numerical costs (remember that finding the leftmost
eigenvalue of the Hessian is the same as finding the global minimizer of the
associated Rayleigh quotient). If we now turn to the computation of $s_k$ in
Step~2, Algorithm~\ref{algo} then computes such a step by attempting to
minimize the model
\beqn{model-p}
T_p(x_k,s) + \frac{\sigma_k}{(p+1)!} \|s\|^{p+1},
\eeqn
as has already been proposed before for
general $p$  \cite{BirgGardMartSantToin17,CartGoulToin17e}.
Moreover, the failure of \req{term-q} in Step~1 is enough, when $q \leq 2$,
to guarantee the existence of nonzero global minimizers of $T_p(x_k,s)$ and
$m_k(s)$, and thus to ensure that a nonzero $s_k$ is possible. The approximate model
minimization is stopped as soon as \req{mterm1} or \req{mterm2} holds, the
latter then reducing to checking that
\beqn{phi-2-m}
\phi_{m_k,2}^\delta(x_k)
= \max\Bigg[
      0,
      - \globmin_{\|d\|\leq \delta}
      \Big(\nabla_s^1m_k(s_k)^Td + \half d^T\nabla_s^2m_k(s_k)d\Big)
      \Bigg]
\leq \frac{\theta\|s_k\|^{p-1}}{(p-1)!}\,\chi_2(\delta)
\eeqn
for some $\delta \in (0,1]$. For each potential $s_k$, finding $\delta \in
  (0,1]$ requires solving (possibly approximately) 
\[
-\globmin_{\|d\|\leq \delta}
\left(\nabla_s^1m_k(s_k)^Td + \half d^T\nabla_s^2m_k(s_k)d\right)
\leq \frac{\theta\|s_k\|^{p-1}}{(p-1)!}\,\chi_2(\delta).
\]
While this could be acceptable without
affecting the overall evaluation complexity of the algorithm, a simpler
alternative is available for $q=2$. We may consider terminating the model
minimization when either \req{mterm1} holds, or
\beqn{termq-2-m}
0>\globmin_{\|d\|\leq 1}
      \Big(\nabla_s^1m_k(s_k)^Td + \half d^T\nabla_s^2m_k(s_k)d\Big)
\geq -\frac{\theta\|s_k\|^{p-1}}{(p-1)!}\,\chi_2(1)
= -\frac{3\theta\|s_k\|^{p-1}}{2(p-1)!}.
\eeqn
The inequality is guaranteed to hold when $s_k$ is close enough
to $s_k^*$, a global minimizer of the model $m_k(s)$, since then
$\nabla_s^1 m_k(s_k^*) = 0$ and $\nabla_s^2 m_k(s_k^*)$ is positive
semi definite, and then $d = 0$ provides the global minimizer
of the second-order Taylor model of $m_k(s)$ around $s_k$.
Verifying \req{termq-2-m} only requires at most one trust-region calculation
for each potential step and ensures \req{phi-2-m} with $\delta=1$, making the
choice $\delta_k=1$ acceptable. The cost this technique is comparable to that
that proposed in \cite{CartGoulToin17e} where an eigenvalue computation is
required for each potential step. Combining these observations,
Algorithm~\ref{algo} then becomes Algorithm~\ref{algo2}.

\algo{algo2}{AR$p$ for $\epsilon$-approximate second-order-necessary minimizers}
{
\vspace*{-0.3 cm}
\begin{description}
\item[Step 0: Initialization.]
  An initial point $x_0\in\calF$  and an initial regularization parameter $\sigma_0>0$
  are given, as well as an accuracy level  $\epsilon \in (0,1)$.  The
  constants $\varpi$, $\theta$, $\eta_1$, $\eta_2$, $\gamma_1$, $\gamma_2$,
  $\gamma_3$ and $\sigma_{\min}$ are also given and satisfy \req{eta-gamma2}.
  Compute $f(x_0)$ and set $k=0$.

\item[Step 1: Test for termination. ]
  Evaluate $\{\nabla^i_x f(x_k)\}_{i=1}^2$.
  If \req{term-q} holds with $\phi_{f,2}^1(x_k)$ given by \req{phi-2} and
  $\delta_{k-1} = 1$, terminate with the approximate solution
  $x_\epsilon=x_k$. Otherwise compute $\{\nabla^i_x f(x_k)\}_{i=3}^p$. 

\item[Step 2: Step calculation. ] 
  Compute a step $s_k \neq 0$  by
  approximately minimizing the model \req{model-p} in the sense that
  \req{descent2} holds and
  \vspace*{-3mm}
  \[
  \|s_k\| \geq \varpi \epsilon^{\frac{1}{p-2+\beta}}
  \tim{ or \req{termq-2-m} holds.}
  \vspace*{-3mm}
  \]

\item[Step 3: Acceptance of the trial point. ]
  Compute $f(x_k+s_k)$ and define $\rho_k$ as in \req{rhokdef2}.
  If $\rho_k \geq \eta_1$, then define
  $x_{k+1} = x_k + s_k$; otherwise define $x_{k+1} = x_k$.

\item[Step 4: Regularization parameter update. ]
  Compute $\sigma_{k+1}$ as in \req{sigupdate2}.
  Increment $k$ by one and go to Step~1 if $\rho_k\geq \eta_1$, or to Step~2 otherwise.
\end{description}
}

\noindent
If $p=q=2$, computing $s_k$ in Step~2 amounts to approximately minimizing the now
well-known cubic model of
\cite{Grie81,NestPoly06,WeisDeufErdm07,CartGoulToin11d}. In
addition, if $s_k$ is the exact global minimizer of this model, the above
argument shows that \req{termq-2-m} automatically holds at $s_k$ and checking
this inequality by solving a trust-region subproblem is thus unnecessary. The only
difference between our proposed algorithm and the more usual cubic
regularization (ARC) method with exact global minimization is that the latter
would check \req{strong-2} for termination, while the
algorithm presented here would instead check \req{phi-2} with $\delta_{k-1}=1$ by
solving a trust-region subproblem. As observed above, both techniques have
comparable numerical cost.

The bound \req{final-bound-Lip} then ensures that Algorithm~\ref{algo2}
terminates in at most $O\Big(\epsilon^{-\frac{p+1}{p-1}}\Big)$ evaluations of
$f$, its gradient and Hessian.  This algorithm thus
shares\footnote{For a marginally weaker (see footnote \ref{gap} and
Theorem~\ref{analog-3.7}) but still necessary and, in our view, more sensible
approximate optimality condition.}  the upper complexity bounds stated in
\cite{CartGoulToin17e} for general $p$ with different values of $\epsilon$ fpr
fisrt- and second-order, and in \cite{NestPoly06,CartGoulToin11d} for $p=2$.

\numsection{A matching lower bound on the evaluation complexity for the
            Lipschitz continuous case}
\label{lowerbound-s}

We now intend to show that the upper bound on evaluation complexity of
Theorem~\ref{upper-theorem} is tight in terms of the order given for
unconstrained and a broad class of constrained problems with Lipschitz
continuous $p$-th derivative (i.e. $\beta = 1$\footnote{A example of slow
convergence for general $\beta$ and $p>1+\beta$ is provided in
\cite{CartGoulToin17d}.}). This objective is attained by defining a variant
of the high-degree Hermite interpolation technique developed in
\cite{CartGoulToin17c}, and then using this technique to build, for any number
$p$ of available derivatives of the objective function and any optimality
order $q$, an unconstrained univariate example of suitably slow convergence
(i.e.\ for which the order in $\epsilon$ given by \req{final-bound-Lip} is
achieved). This example is then embedded in higher dimensions to provide
general lower bounds.

\subsection{High-degree univariate Hermite interpolation}\label{Hermite-s}

We start by investigating some useful properties of Hermite interpolation.
Let us assume that we wish to construct a univariate Hermite interpolant
$\pi$ of degree $2(p+1)$ of the form
\beqn{pi-def}
\pi(\tau) = \sum_{i=0}^{2p+1} c_i \, \tau^i
\eeqn
on the interval $[0,s]$ satisfying the $2(p+1)$ conditions
\beqn{interp-conds}
\pi^{(i)}(0) = f^{(i)}_0, \ms \pi^{(i)}(s)= f^{(i)}_1
\tim{for} i\in \iibe{0}{p},
\eeqn
where $f^{(i)}_0$ and $f^{(i)}_1$ are given. The values of the coefficients
$c_0,\ldots , c_p$ may then be obtained by
\[
c_i = \frac{f^{(i)}_0}{i!} \tim{for} i \in \iibe{0}{p}
\]
while the remaining ones satisfy the linear system
\beqn{interpsys-1}
\matr{ccccc}{a_{0,0} s^{p+1}  & a_{0,1}s^{p+2} & \cdots & a_{0,p-1}s^{2p}   & a_{1,p}s^{2p+1} \\
             a_{1,0}s^p      & a_{2,2}s^{p+1} & \cdots & a_{2,p-1} s^{2p-1} & a_{2,p}s^{2p}\\
             \vdots         &    \vdots    & \ddots &     \vdots       & \vdots \\
             a_{p,0} s       &  a_{p,1}s^2   & \cdots & a_{p,p-1}s^p      & a_{p,p}s^{p+1} }
\cvect{c_{p+1} \\ c_{p+2} \\ \vdots \\ c_{2p+1} }
= \cvect{
  f_1^{(0)} - T_p^{(0)}(0,s) \\ f^{(1)}_1 - T_p^{(1)}(0,s)\\ \vdots \\ f^{(p)}_1 - T_p^{(p)}(0,s) )
  }
\eeqn
where
\[
T_p(0,s) = \sum_{i=0}^p \frac{f_0^{(i)}}{i!}\, s^i
\tim{ and }
a_{i,j}
= \frac{(p+j+1)!}{(p+j+1-i)!}
\ms (i,j=0, \ldots, p).
\]
Observe that \req{interpsys-1} can be rewritten as
\[
\matr{ccccc}{ s^p & &  & & \\
  & s^{p-1}& & & \\
  & \ddots & & \\
  & & & & 1} A_p
\matr{ccccc}{ s &  & & & \\
  & s^2 & & & \\
  & & \ddots & & \\
  & & & & s^{p+1}}
 \cvect{c_{p+1} \\ c_{p+2}  \\ \vdots \\ c_{2p+1} }
= \cvect{
  f_1^{(0)} - T_p^{(0)}(0,s) \\
  f_1^{(1)} - T_p^{(1)}(0,s) \\
  \vdots \\
  f^{(p)}_1 - T_p^{(p)}(0,s) )
}
\]
with $A_p$ is the matrix whose $(i,j)$-th entry is $a_{i,j}$, which only
depends on $p$. It was show in \cite[Appendix]{CartGoulToin17c} that $A_p$ is
nonsingular.  Therefore
\[
 \cvect{c_{p+1}\,s \\ c_{p+2}\,s^2 \\\vdots \\ c_{2p+1}\, s^{p+1} }
= A_p^{-1} \cvect{
  \frac{1}{s^p}[f_1^{(0)} - T_p^{(0)}(0,s)] \\
  \frac{1}{s^{p-1}}[f_1^{(1)} - T_p^{(1)}(0,s)] \\
  \vdots \\
  f^{(p)}_1 - T_p^{(p)}(0,s) }.
\]
We therefore deduce that, for any $\tau \in [0,s]$ ,
\[
\begin{array}{lcl}
|\pi^{(p+1)}(\tau)|
& = & \left|\bigsum_{i=0}^p \bigfrac{(p+1+i)!}{i!}c_{p+1+i}\, \tau^i \right|\\
& \leq & \bigsum_{i=0}^p \bigfrac{(p+1+i)!}{i!}\left(| c_{p+1+i}|\, s^{i+1} \right) s^{-1}\\
& \leq & \bigfrac{(p+1)(2p+1)!}{p!} \|A_p^{-1}\|_\infty
         \bigmax_{j=0,\ldots,p}\left|\frac{f_1^{(j)} - T_p^{(j)}(0,s)}{s^{p-j+1}}\right|.
\end{array}
\]
The mean-value theorem then implies that, for any
$0\leq \tau_2 \leq \tau_1 \leq s$ and some
$\xi \in [\tau_2,\tau_1]\subseteq [0,s]$,
\beqn{pi-p-aH}
\begin{array}{lcl}
\bigfrac{|\pi^{(p)}(\tau_1)-\pi^{(p)}(\tau_2)|}{|\tau_1-\tau_2|}
&   =  & |\pi^{(p+1)}(\xi) | \\
& \leq & \bigmax_{\tau \in [0,s]}|\pi^{(p+1)}(\tau) | \\
& \leq & \bigfrac{(p+1)(2p+1)!}{p!} \|A_p^{-1}\|_\infty
         \bigmax_{j=0,\ldots,p}\left|\frac{f_1^{(j)} - T_p^{(j)}(0,s)}{s^{p-j+1}}\right|.
\end{array}
\eeqn
This development thus leads us to the following conclusion.

\lthm{interp-Holder}{Suppose that $\{f_\ell^{(j)}\}$ are given for $\ell \in \{1,2\}$ and $j \in
  \iibe{0}{p}$. Suppose also that there exists a constant
  $\kappa_f \geq 0$ such that, for all $j \in \iibe{0,}{p}$,
  \beqn{data-Holder}
  |f_1^{(j)} - T_p^{(j)}(0,s)| \leq \kappa_f \, s^{p-j+1}.
  \eeqn
  Then the Hermite interpolation polynomial $\pi(\tau)$ on $[0,s]$
  given by \req{pi-def} and satisfying \req{interp-conds} admits a
  Lipschitz continuous $p$-th derivative on $[0,s]$, with Lipschitz constant
  given by
  \[
  L_p \eqdef \bigfrac{(p+1)(2p+1)!}{p!} \|A_p^{-1}\|_\infty \kappa_f,
  \]
  which only depends on $p$ and $\kappa_f$.
}

\proof{Directly results from \req{pi-p-aH} and \req{data-Holder}.
}

\noindent
Observe that \req{data-Holder} is identical to \req{resder} when $\beta=1$ and
$n=1$. This means that the conditions of Theorem~\ref{interp-Holder} automatically
hold if the interpolation data $\{f_i^{(j)}\}$ is itself extracted
from a function having a Lipschitz continuous $p$-th derivative.

Applying the above results to several interpolation intervals then yields the
existence of a smooth Hermite interpolant.

\lthm{global-pi-th}{
  Suppose that, for some integer $k_e >0$ and $p>0$, the data
  $\{f_k^{(j)}\}$ and $\{x_k\}$ is given for $k \in \iibe{0}{k_e}$ and $j \in
  \iibe{0}{p}$.  Suppose also that $s_k= x_{k+1}-x_k \in (0,\kappa_s]$ for $k \in
\iibe{0}{k_e}$ and some $\kappa_s >0$, and that, for some
  constant $\kappa_f \geq 0$ and $k \in \iibe{0}{k_e-1}$,
  \beqn{data-aH}
  |f_{k+1}^{(j)} - T_{k,p}^{(j)}(x_k,s_k)| \leq \kappa_f \, s_k^{p-j+1}.
  \eeqn
  where $T_{k,p}(x_k,s) = \sum_{i=0}^p f_k^{(i)}\, s^i/i!$.  Then there exists a
  $p$ times continuously differentiable function $f$ from $\Re$ to $\Re$ with
  Lipschitz continuous $p$-th derivative such that, for $k \in
  \iibe{0}{k_e}$, 
  \[
  f^{(j)}(x_k) = f_k^{(j)} \tim{for} j\in\iibe{0}{p}.
  \]
  Moreover, the range of $f$ only depends on $p$, $\kappa_f$, $\max_k f_k^{(0)}$ and
  $\min_k f_k^{(0)}$.
}

\proof{
  We first use Theorem~\ref{interp-Holder} to define a Hermite interpolant
  $\pi_k(s)$ of the form \req{pi-def} on each interval $[x_k,
    x_{k+1}]=[x_k,x_k+s_k]$ ($k \in \iibe{0}{k_e}$)
  using $f_0^{(j)} = f_k^{(j)}$ and $f_1^{(j)}= f_{k+1}^{(j)}$ for $j\in\iibe{0}{p}$,
  and then set
  \[
  f(x_k+s) = \pi_k(s)
  \]
  for any $s \in [0,s_k]$.  We may then smoothly prolongate $f$ for $x\in
  \Re$ by defining two additional interpolation intervals
  $[x_{-1},x_0]=[-s_{-1},0]$ and  $[x_{k_e},x_{k_e}+s_{k_e}]$   with end conditions
  \[
  f_{-1} = f_0^{(0)},
  \ms
  f_{k_e+1} = f_{k_e}^{(0)}
  \tim{and}
  f_{-1}^{(j)} = f_{k_e+1}^{(j)} = 0
  \tim{for} j \in \ii{p},
  \]
  and where $s_{-1}$ and $s_{k_e}$ are chosen sufficiently large to ensure that
  \req{data-aH} also holds on intervals -1 and $k_e$.  We next set
  \[
  f(x) = \left\{  \begin{array}{ll}
         f_0^{(0)} & \tim{for } x \leq x_{-1},\\
         \pi_k(x-x_k) & \tim{for } x\in[x_k,x_{k+1}]
                                 \tim{and} k\in \iibe{-1}{n},\\
         f_{k_e}^{(0)} & \tim{for } x \geq x_{k_e}+s_{k_e}.\\
         \end{array}\right.
   \]
}    

\subsection{Slow convergence to ($\epsilon$,$\delta$)-approximate $q$-th-order-necessary
  minimizers}
\label{examples-s}

We now consider an unconstrained univariate instance of problem
\req{problem}. Our aim is first to show that, for each choice of $p\geq 1$ and
$q\in\ii{p}$, there exists an objective function $f$ for problem \req{problem}
with $f \in C^{p,1}(\Re)$ (i.e. $\beta=1$) such that obtaining an
$(\epsilon,\delta)$-approximate $q$-th-order-necessary minimizer may require 
at least
\[
\epsilon^{-\frac{p+1}{p-q+1}}
\]
evaluations of the objective function and its derivatives using
Algorithm~\ref{algo}, matching, in order of $\epsilon\in (0,1]$, the upper bound
\req{final-bound-Lip}. Our development follows the broad outline of
\cite{CartGoulToin18a} but extends it to approximate minimizers of arbitrary
order. Given a model degree $p \geq 1$ and an optimality order $q\in \ii{p}$,
we first define the sequences $\{f_k^{(j)}\}$  for $j \in
\iibe{0}{p}$ and $k\in\iibe{0}{k_\epsilon}$ with 
\beqn{keps}
k_\epsilon = \left\lceil \epsilon^{-\frac{p+1}{p-q+1}}\right\rceil
\eeqn
by
\beqn{reg-omegak}
\omega_k= \epsilon \, \frac{k_\epsilon-k}{k_\epsilon}.
\eeqn
as well as
\beqn{fnotq-def}
f_k^{(j} = 0
\tim{for} j \in \ii{q-1} \cup \iibe{q+1}{p}
\eeqn
and
\beqn{fq-def}
f_k^{(q)} = -(\epsilon+\omega_k)\,q!\,\chi_q(1) < 0.
\eeqn
Thus
\beqn{4.10a}
T_p(x_k,s) =
\sum_{j = 0}^p \frac{f_k^{(j)}}{j!} s^j
= f_k^{(0)} -(\epsilon+\omega_k) \chi_q(1) s^q
\eeqn
and, assuming $\delta_{k-1}=1$ for all $k$ (we verify below that this is
acceptable),
\beqn{phiex}
\phi_{f,q}^{\delta_{k-1}}(x_k) = (\epsilon+\omega_k) \chi_q(\delta_{k-1})
\eeqn
We also set $\sigma_k = p!$ for all $k\in \iibe{0}{k_\epsilon}$ (we again
verify below that is acceptable).  Note that
\beqn{4.9}
\omega_k \in (0,\epsilon]
\tim{ and }
\phi_{f,q}^{\delta_{k-1}}> \epsilon \chi_q(\delta_{k-1})
\tim{for}
k \in \iibe{0}{k_\epsilon-1},
\eeqn
(and \req{term-q} fails at $x_k$), while
\beqn{4.10}
\omega_{k_\epsilon} =  0
\tim{ and }
\phi_{f,q}^{\delta_{k-1}}(x_{k_\epsilon}) = \epsilon \chi_q(\delta_{k-1})
\eeqn
(and \req{term-q} holds at $x_{k_\epsilon}$). It is easy to verify using
\req{4.10a} that the model \req{model} is then globally minimized for
\beqn{sstar-def}
s_k = \left[ \frac{|f^{(q)}_k|}{(q-1)!}\right]^{\frac{1}{p-q+1}}
= [q(\epsilon+\omega_k)\chi_q(1)]^{\frac{1}{p-q+1}}
> \epsilon^{\frac{1}{p-q+1}}
\ms
(k\in\iibe{0}{k_\epsilon}).
\eeqn
Hence this step satisfies 
\req{mterm1} if we choose $\varpi=1$.  Because of this fact, we are free to
choose $\delta_k$ arbitrarily in $(0,1]$ and we choose $\delta_k = 1$. Thus,
  provided we make the choice $\delta_{-1}=1$ ensuring \req{phiex} for $k=0$ ,
  the value $\delta_k=1$ is admissible for all $k$. The step \req{sstar-def}
  yields that
\beqn{4.12}
\begin{array}{lcl}
  m_k(s_k)
  &=& f_k^{(0)} -
  (\epsilon+\omega_k)\chi_q(\delta_k)[q(\epsilon+\omega_k)\chi_q(\delta_k)]^{\frac{q}{p-q+1}}
  + \frac{1}{p+1}[q(\epsilon+\omega_k)\chi_q(\delta_k)]^{\frac{p+1}{p-q+1}}\\*[2ex]
  &=& f_k^{(0)} - \zeta(q,p)[q(\epsilon+\omega_k)\chi_q(\delta_k)]^{\frac{p+1}{p-q+1}}
\end{array}
\eeqn
where
\beqn{4.13}
\zeta(q,p) \eqdef \frac{p-q+1}{q(p+1)} \in (0,1).
\eeqn
Thus $m_k(s_k)<m_k(0)$ and \req{descent2} holds.
We then define
\beqn{4.14}
f_0^{(0)} = 2[2q\chi_q(1)]^{\frac{p+1}{p-q+1}}
\tim{ and }
f_{k+1}^{(0)} =
f_k^{(0)}-\zeta(q,p)[q(\epsilon+\omega_k)\chi_q(\delta_k)]^{\frac{p+1}{p-q+1}},
\eeqn
which provides the identity
\beqn{4.15}
m_k(s_k) = f_{k+1}^{(0)}
\eeqn
(ensuring that iteration $k$ is successful because $\rho_k=1$ in
\req{rhokdef2} and thus that our choice of a constant $\sigma_k$ is
acceptable). In addition, using \req{4.14}, \req{4.9}, \req{4.13}, the
equality $\delta_k=1$ and the inequality $k_\epsilon\leq
1+\epsilon^{-\frac{p+1}{p-q+1}}$ from \req{keps} gives that, for
$k \in \iibe{0}{k_\epsilon}$,
\[
\begin{array}{lcl}
f_0^{(0)} \geq f_k^{(0)}
& \geq & f_0^{(0)}-k\zeta(q,p)[2q\epsilon\chi_q(\delta_k)]^{\frac{p+1}{p-q+1}}\\*[2ex]
& \geq & f_0^{(0)}-k_\epsilon\epsilon^{\frac{p+1}{p-q+1}}[2q\chi_q(1)]^{\frac{p+1}{p-q+1}}\\*[2ex]
& \geq & f_0^{(0)}-\Big(1+\epsilon^{\frac{p+1}{p-q+1}}\Big)[2q\chi_q(1)]^{\frac{p+1}{p-q+1}}\\*[2ex]
& \geq & f_0^{(0)}-2[2q\chi_q(1)]^{\frac{p+1}{p-q+1}},
\end{array}
\]
and hence that
\beqn{4.16}
f_k^{(0)} \in \left[ 0, 2[2q\chi_q(1)]^{\frac{p+1}{p-q+1}} \right]
\tim{ for }
k \in \iibe{0}{k_\epsilon}.
\eeqn
We also set
\[
\delta_{-1} = 1,
\ms
x_0 = 0
\tim{ and }
x_k = \sum_{i=0}^{k-1} s_i.
\]
Then \req{4.15} and \req{model} give that
\beqn{conds-0}
|f^{(0)}_{k+1}- T_p(x_k,s_k)| = \frac{1}{p+1}|s_k|^{p+1}.
\eeqn
Now note that, using \req{4.10a} and the first equality in \req{sstar-def},
\[
T_p^{(j)}(x_k,s_k)
= \frac{f_k^{(q)}}{(q-j)!} \,s_k^{q-j}\,\delta_{[j\leq q]}
= -\bigfrac{(q-1)!}{(q-j)!}\,s_k^{p-j+1}\,\delta_{[j\leq q]}
\]
where $\delta_{[\cdot]}$ is the standard indicator function.
We may now verify that, for $j\in\ii{q-1}$, 
\beqn{conds-noq1}
| f^{(j)}_{k+1}- T^{(j)}_p(x_k,s_k)|
= |0-T^{(j)}_p(x_k,s_k)|
\leq \left|\frac{(q-1)!}{(q-j)!}\right|\,|s_k|^{p-j+1}
\leq (q-1)! \,|s_k|^{p-j+1},
\eeqn
while, for $j=q$, we have that
\beqn{conds-q}
| f^{(q)}_{k+1}- T^{(q)}_p(x_k,s_k)|
= | -(q-1)!\,s_k^{p-q+1} + (q-1)!\,s_k^{p-q+1}|
= 0
\eeqn
and, for $j \in \iibe{q+1}{p}$,
\beqn{conds-noq2}
| f^{(j)}_{k+1}- T^{(j)}_p(x_k,s_k)|
= |0 - 0|
= 0.
\eeqn
Combining \req{conds-0}, \req{conds-noq1}, \req{conds-q} and \req{conds-noq2},
we deduce that \req{data-aH} holds with $\kappa_f = (q-1)!$. We
may thus apply Theorem~\ref{global-pi-th} with $\beta=1$, $\kappa_f = (q-1)!$
and $\kappa_s=1$, and deduce the existence of a $p$ times continuously
differentiable function $f$ from $\Re$ to $\Re$ with Lipschitz continuous
derivatives of order $0$ to $p$ which interpolates the $\{f_k^{(j)}\}$ at
$\{x_k\}$ for $k \in \iibe{0}{n}$ and $j \in \iibe{0}{p}$.  Moreover,
\req{4.16} and Theorem~\ref{global-pi-th} imply that the range of $f$
only depends on $p$ and $q$. In addition, \req{4.15} ensures that every
iteration is successful and thus, because of \req{sigupdate2}, that the value
$\sigma_k=p!$ may be used at all iterations.

This argument allows us to state the following lower bound on the complexity
of the regularization algorithm using a $p$-th degree model.

\llem{lower-theorem}{
  Given any  $p \in \Na_0$ and $q\in\ii{p}$, there exists a $p$
  times continuoulsy differentiable function $f$ from $\Re$ to $\Re$ with
  range only depending on $p$ and $q$ and Lipschitz continuous $p$-th
  derivative such that, when the regularization algorithm with $p$-th degree
  model (Algorithm~\ref{algo}) is applied to minimize $f$ without constraints,
  it takes exactly
  \[
  k_\epsilon = \left\lceil \epsilon^{-\frac{p+1}{p-q+1}}\right\rceil
  \]
  iterations (and evaluations of the objective function and its derivatives)
  to find an ($\epsilon$,$\delta)$-approximate $q$-th-order-necessary minimizer.
}

\noindent
This implies the following important consequence for higher dimensional problems.

\lthm{multivariate-lower-bound}{
  Given any $n \in \Na_0$, $p \in \Na_0$ and $q\in\ii{p}$,
  there exists a $p$ times continuoulsy differentiable function $f$ from
  $\Re^n$ to $\Re$ with range only depending on $p$ and $q$ and Lipschitz
  continuous $p$-th derivative tensor such that, when  the regularization
  algorithm with $p$-th degree model (Algorithm~\ref{algo}) is applied to
  minimize $f$ without constraints, it takes exactly 
  \beqn{keps-lower}
  k_\epsilon = \left\lceil \epsilon^{-\frac{p+1}{p-q+1}}\right\rceil
  \eeqn
  iterations (and evaluations of the objective function and its derivatives)
  to find an  $(\epsilon,\delta)$-approximate $q$-th-order-necessary minimizer.
  Furthermore, the same conclusion holds if the optimization
  problem under consideration involves constraints provided the feasible set
  $\calF$ contains a ray.
}

\proof{
The first conclusion directly follows from Lemma~\ref{lower-theorem} since it is
always possible to include the unimodal example as an independent component of
a multivariate one.

The second conclusion follows from the observation that our univariate example
of slow convergence is only defined on $\Re^+$ (even if Theorem~\ref{global-pi-th}
provides an extension to the complete real line). As a consequence, it may be
used on any feasible ray.
}

We now make a few observations.
\begin{enumerate}
\item Theorem~\ref{multivariate-lower-bound} generalizes to arbitrary $q$ the
  bound obtained in \cite{CarmDuchHindSidf17a} for the case $q=1$ and also
  shows that, at variance with the result derived in this reference, the
  generalized bound applies for arbitrary problem's dimension, but depends on
  $\epsilon$, $p$ and $q$.
  
\item For simplicity, we have chosen, in the above example, to minimize the
  model $m_k(s)$ globally at every iteration, but we might consider
  other pairs $(s_k,\delta_k)$. A similar example of slow convergence may in
  fact be constructed along the lines used above\footnote{At the price of
  possibly larger constants.}  for any sequence of
  acceptable\footnote{Remember that $\delta=1$ is always possible for
  $q=1$. It thus unsurprising that no such condition appears 
  in \cite{CarmDuchHindSidf17a}.} model reducing steps and
  associated optimality radii (in the sense of Lemma~\ref{step-ok-l}),
  provided the optimality radii remain bounded away from zero.  This means
  that our example of slow convergence applies not only to
  Algorithm~\ref{algo} but also to a much broader class of minimization
  methods.  Moreover, it is also possible to weaken the constraints on the
  step further by relaxing \req{4.15} and only insisting on acceptable
  decrease of the objective function value in Step~3 of the algorithm.

  In \cite{CarmDuchHindSidf17a}, the authors derive their upper bound for
  $q=1$ for the general class of ``zero-preserving'' algorithms, which are
  algorithms that ``never explore (from $x_k$) coordinates which appear not to
  affect the function'', that is directions $d$ along which $T_p(x_k,\cdot)$
  is constant. This property is obviously shared by Algorithm~\ref{algo} because it
  attempts to reduce the Taylors' expansion of $f$ around the current iterate
  (the presence of the isotropic regularization term is irrelevant for this).

\item Our example does not apply, for instance, to a linesearch method using
  global univariate minimization in a direction of search computed from the
  Taylor's expansion of $f$, which is another zero-preserving method. Note
  however that this method, just as every other linesearch method (including
  possibly randomized coordinate searches), is bound to fail when attempting
  to compute approximate minimizers of order beyond three, because the
  Taylor's expansion at a non-optimal point then needs no longer decrease
  along lines.  This is demonstrated by the following old example
  \cite{Hanc17,Pean84}.  Let 
  \[
  f(x_1,x_2) = (\half x_1^2 - x_2)(x_1^2 - x_2).
  \]
  Then $f(0,0)= 0$ and the origin is not a minimizer since $f$
  decreases along the arc $x_2=\threequarters x_1$. Yet the origin is the
  global minimizer along every line passing through the origin, preventing
  any linesearch method to progress away from $(0,0)$.
\end{enumerate}

\noindent
Let us now consider an alternative unconstrained minimization method which
would attempt to reduce the \emph{unregularized} model (that is \req{model}
with $\sigma_k = 0$) in order to find an unconstrained first-order minimizer.
It is easy to see that if one chooses 
\[
f^{(1)}_k = -(\epsilon+\omega_k),
\ms
f^{(i)}_k = 0 \tim{for} i\in\iibe{2}{p-1}
\tim{and}
f^{(p)}_k = p!,
\]
the same reasoning as above yields that the largest obtainable decrease with
this model occurs at
\[
s_k = \left(\frac{\epsilon+\omega_k}{p}\right)^{\frac{1}{p-1}}
\]
and is given by
\[
f^{(0)}_k- m_k(s_k)= (p-1)\left(\frac{\epsilon+\omega_k}{p}\right)^{\frac{p}{p-1}}.
\]
This then implies that at least a multiple of $\epsilon^{-\frac{p}{p-1}}$
evaluations may be needed to find approximate first-order-necessary
minimizers, which is worse than the bound in $\epsilon^{-\frac{p+1}{p}}$
holding for the regularized algorithm. This is consistent with the known lower
$O(\epsilon^{-2})$ bound for first-order points that holds for the
(unregularized) Newton method (and hence the trust-region method), both of
which use $p=2$. Adding the regularization term thus not only provides a
mechanism to limit the stepsize and make the step well-defined when
$T_p(x_k,s)$ is unbounded below, but also amounts to increasing the 'useful
degree' of the model by one, improving the worst-case complexity bound.

Summing up the above discussion, we conclude that an example of slow
convergence requiring at least \req{keps-lower} evaluations can be built for
any method whose steps decrease the regularized ($\sigma_k\geq \sigma_{\min}$)
or unregularized ($\sigma_k=0$) model \req{model} and whose approximate local
optimality can be measured by \req{mterm2} for some constant $\theta$ and
$\delta_k=1$ (which we can always enforce by adapting $\varpi$ and 
\req{fnotq-def}). For orders up to two, this includes most variants of
steepest-descent and Newton's methods including those
globalized with regularization, trust-region, a linesearch or a mixture of these
(see \cite{CartGoulToin18a} for a discussion). General linesearch methods are
excluded for high-order optimization as they may fail to converge to
approximate minimizers of order four and beyond.

Finally, one may wonder at what would happen if, for the
interpolation data \req{fnotq-def}-\req{fq-def}, the model
\[
m_k(s) = T_p(x_k,s)+ \frac{\sigma_k}{m!}|s|^m
\]
were used for some $m > p+1$, resulting in a shorter step. The global model
minimizer would then occur at $s = [q(\epsilon+\omega_k)\chi_q(1)]^{1/(m-1)}$
and give an optimal model decrease equal to
$[q(\epsilon+\omega_k)\chi_q(1)]^{m/(m-1)}(m-q)/m$.  However,
\req{data-aH} would then fail for $j=0$ and the argument leading to an example
of slow convergence would break down.

\numsection{Summary, further comments and open questions}\label{concl-s}

For any optimality order $q \geq 1$, we have provided the concept of an
$(\epsilon,\delta)$-approximate $q$-th-order-necessary minimizer for the very
general set-constrained problem \req{problem}.  We have then proposed a
conceptual regularization algorithm to find such approximate minimizers and
have shown that, if $\nabla_x^p f$ is $\beta$-H\"{o}lder continuous, this
algorithm requires at most $O(\epsilon^{-\frac{p+\beta}{p-q+\beta}})$
evaluations of the objective function and its $p$ first derivatives to
terminate.  When $\nabla_x^p f$ is Lipschitz continuous, we have used an
unconstrained univariate version of the problem to show that this bound is
sharp in terms of the order in $\epsilon$ for any feasible set containing a
ray and any problem dimension.

In view of the results in \cite{CartGoulToin16,NestGrap16}, one may wonder at
what would happen if the regularization power (i.e.\ the power of $\|s\|$ used
in the last term of the model \req{model}) is allowed to differ from
$p+\beta$. The theory presented above must then be re-examined and the crucial
point is whether a global upper bound $\sigma_{\max}$ on the regularization
parameter can still be ensured as in Lemma~\ref{sigmaupper-lemma}.  One easily
verifies that this is the case for regularization powers $r\in
(p,p+\beta]$. Arguments parallel to those presented above then yield an upper 
bound of $O(\epsilon^{-\frac{r}{r-q}})$ evaluations\footnote{We may even relax
\req{mterm2} slightly by replacing $\|s_k\|^{p-q+\beta}$ by
$\|s_k\|^{r-q}$.}, recovering the bound given in Section~3.3 of
\cite{CartGoulToin16} for $q=1$. The situation is however more complicated
(and beyond the scope of the present paper) for $r > p+\beta$ and the
determination of a suitable general complexity upper bound for this latter
case has not been formalized at this stage, but the analysis for $q=1$
discussed in Section~3.2 of \cite{CartGoulToin16} suggests that an improvement
of the bound for larger $r$ is unlikely.

Although the results presented essentially solve the question of determining
the optimal evaluation complexity for unconstrained problems and problems with
general inexpensive constraints, some interesting issues remain open at this
stage. A first such issue is whether an example of slow convergence for all
$\epsilon \in(0,1)$ can be found for feasible domains not containing a ray. A
second is to extend the general complexity theory for problems whose
constraints are not inexpensive: the discussion in \cite{CartGoulToin17b}
indicates that this is a challenging research area.

{\footnotesize


}


\newcommand{\resetcounters}{\setcounter{equation}{0} \setcounter{figure}{0}
 \setcounter{table}{0}}
\newcommand{\newappendixname}{A}
\newcommand{\newappendix}[1]{\renewcommand{\newappendixname}{{#1}}
 \section*{Appendix \newappendixname}
 \resetcounters
 \renewcommand{\theequation}{\newappendixname.\arabic{equation}}
 }

\renewcommand\appendix{%
  \setcounter{section}{0}%
  \setcounter{subsection}{0}%
  \renewcommand\thesection{\Alph{section}}}

\renewcommand{\thesection}{Appendix \Alph{section}}
\setcounter{section}{0}%
\setcounter{subsection}{0}%


\newappendix{A}
\renewcommand{\theequation}{A.\arabic{equation}}
\renewcommand{\thesection}{A.\arabic{section}}
\renewcommand{\thesubsection}{A.\arabic{subsection}}

\section{Proof of Lemmas in Section \ref{arp-s}}

\noindent

{\bf Proof of Lemma \ref{taylor-bounds-lemma}.}\label{A-Lemma-proof}
We first establish the identity
\beqn{intid}
I_{k-1,\beta} \eqdef
\bigint_0^1 \xi^{\beta} (1-\xi)^{k-1} \, d\xi = \frac{(k-1)!}{(k+\beta)!},
\tim{where} (k+\beta)! \eqdef \prod_{i=1}^k (i+\beta).
\eeqn
To see this, integrating by parts, we have that
\[
I_{k-1,\beta} = \left[ -(k-1) \xi^{\beta} (1-\xi)^{k-2} \right]_0^1
+ \frac{(k-1)}{(1+\beta)} \xi^{1+\beta}(1-\xi)^{k-2} \, d\xi
= \frac{(k-1)}{(1+\beta)} I_{k-2,1+\beta}
\]
and thus, recursively, that
\[
I_{k-1,\beta} = \frac{(k-1)!}{(k-1+\beta)!} I_{0,k-1+\beta}
=  \frac{(k-1)!}{(k-1+\beta)!} \int_0^1 \xi^{k-1+\beta} \, d\xi
=  \frac{(k-1)!}{(k+\beta)!}.
\]
As in \cite{CartGoulToin17c}, consider the Taylor identity
\beqn{unitaylor}
\psi(1) - \tau_k(1) = \frac{1}{(k-1)!}\int_0^1 (1 - \xi)^{k-1}
                           [\psi^{(k)}(\xi ) - \psi^{(k)}(0)]  \, d\xi
\eeqn
involving a given univariate $C^k$ function $\psi(t)$
and its $k$-th order Taylor
approximation
\[
\tau_k(t) = \sum_{i=0}^k \psi^{(i)}(0) \frac{t^i}{i!}
\]
expressed in terms of the value $\psi^{(0)} = \psi$ and
$i$th derivatives $\psi^{(i)}$, $i=1,\ldots,k$.
Then, picking $\psi(t) = f(x+t s)$, for given
$x, s \in \Re^n$, and $k = p$,
the identity \req{unitaylor},
and the relationships $\psi^{(p)}(t) = \nabla_x^p f(x+t s)[s]^p$
and $\tau_p(1) = T_p(x,s)$ give that
\[
f(x+s) -  T_p(x,s) = \frac{1}{(p-1)!}\int_0^1 (1 - \xi)^{k-1}
 \left( \nabla_x^p f(x+\xi s) - \nabla_x^p f(x)\right)[s]^p \, d\xi,
\]
and thus from the definition of the tensor norm \req{Tnorm},
the H\"{o}lder bound \req{tensor-Hol} and the identity
\req{intid} when $k=p$ that
\[
\begin{array}{ll}
  f(x+s) - T_p(x,s) \!\!\!
  & \leq  \bigfrac{1}{(p-1)!}\bigint_0^1  (1 - \xi)^{k-1}
  \left| \left(\nabla_x^p f(x+\xi s) - \nabla_x^p f(x)\right)
  \left[\frac{s}{\|s\|}\right]^p \right|
    \|s\|^p \, d\xi

    \\ & \leq  \bigfrac{1}{(p-1)!}\bigint_0^1  (1 - \xi)^{k-1} \max_{\|v\|=1}
    \left| \left(\nabla_x^p f(x+\xi s) - \nabla_x^p f(x)\right)
  \left[v\right]^p \right|
    \|s\|^p \, d\xi
    \\ & = \bigfrac{1}{(p-1)!}\bigint_0^1  (1 - \xi)^{k-1}
    \| \nabla_x^p f(x+\xi s) - \nabla_x^p f(x))\|_{[p]} d\xi \cdot
    \|s\|^p \,
  \\ &   \leq \bigfrac{1}{(p-1)!}
  \bigint_0^1 \xi^{\beta} (1-\xi)^{p-1} \, d\xi \cdot L \|s\|^{p+\beta}
    = \bigfrac{L}{(p+\beta)!} \|s\|^{p+\beta}
\end{array}
\]
for all $x,s \in \Re^n$, which is the required \req{resf}.

Likewise, for arbitrary unit vectors $v_1,\ldots,v_j$,
choosing $\psi(t) = \nabla_x^j f(x+t s)[v_1,\ldots,v_j]$
  and $k=p-j$,it follows from \req{unitaylor},
the relationships $\psi^{(p-j)}(t) =
 \nabla_x^p f(x+t s)[v_1,\ldots,v_j][s]^{p-j}$
 and $\tau_{p-j}(1) = \nabla_s^j T_p(x,s)$ that
 \beqn{jdiff}
 \begin{array}{ll}
  ( \nabla_x^j  \!\!\!\!\! & f(x+s) -  \nabla_s^j T_p(x,s)) [v_1,\ldots,v_j] \\
& = \bigfrac{1}{(p-j-1)!} \bigint_0^1 (1-\xi)^{p-j-1} \left(
   \nabla_x^p f(x+\xi s) -  \nabla_x^p f(x) \right [v_1,\ldots,v_j][s]^{p-j}
   \, d\xi.
\end{array}
\eeqn
Then picking $v_1,\ldots,v_j$ to maximize the absolute value of
left-hand size of \req{jdiff} and using
 the tensor norm \req{Tnorm}, the H\"{o}lder bound \req{tensor-Hol} and
 the identity \req{intid} when $k=p-j$, we find that
\[
\begin{array}{ll}
\| \nabla_x^j  \!\!\!\!\! & f(x+s) -  \nabla_s^j T_p(x,s) \|_{[j]}
 \\ & \leq \bigfrac{1}{(p-j-1)!} \bigint_0^1 (1-\xi)^{p-j-1}
 \left| ( \nabla_x^p f(x+\xi s) -  \nabla_x^p f(x)) [v_1,\ldots,v_j]
  \left[\bigfrac{s}{\|s\|}\right]^{p-j} \right| \|s\|^{p-j} \, d\xi
 \\ & \leq \bigfrac{1}{(p-j-1)!} \bigint_0^1 (1-\xi)^{p-j-1} \!\!\!
\max_{\|v_1\|=\cdots=\|v_p\|=1}
\left| \left ( \nabla_x^p f(x+\xi s) -  \nabla_x^p f(x) \right )
       [v_1,\ldots,v_p] \right| \|s\|^{p-j} \, d\xi
 \\ & = \bigfrac{1}{(p-j-1)!} \bigint_0^1 (1-\xi)^{p-j-1}
\| \nabla_x^p f(x+\xi s) -  \nabla_x^p f(x) \|_{[p]} \, d\xi \cdot \|s\|^{p-j}
 \\ & \leq \bigfrac{1}{(p-j-1)!}
 \bigint_0^1 \xi^{\beta} (1-\xi)^{p-j-1} \, d\xi \cdot L \|s\|^{p-j+\beta}
= \bigfrac{L}{(p-j+\beta)!} \|s\|^{p-j+\beta}
\end{array}
\]
for all $x,s \in \Re^n$,
which gives \req{resder}.\hfill$\Box$\\

\noindent
{\bf Proof of Lemma \ref{SvsU}.}
The regularization parameter update \req{sigupdate2}
gives that, for each $k$,
\[
\gamma_1\sigma_j \leq \max[\gamma_1\sigma_j,\sigma_{\min}]
 \leq \sigma_{j+1}, \ms j \in \calS_k,
\tim{ and }
\gamma_2\sigma_j \leq \sigma_{j+1}, \ms j \in \calU_k,
\]
where $\calU_k \eqdef \iibe{0}{k}\setminus\calS_k$.
Thus we deduce inductively that
$\sigma_0\gamma_1^{|\calS_k|}\gamma_2^{|\calU_k|}\leq \sigma_{k}$.
We therefore obtain, using \req{sigmax}, that
\[
|\calS_k|\log{\gamma_1}+|\calU_k|\log{\gamma_2}\leq
\log\left(\frac{\sigma_{\max}}{\sigma_0}\right),
\]
which then implies that
\[
|\calU_k|
\leq - |\calS_k|\frac{\log\gamma_1}{\log\gamma_2}
      + \frac{1}{\log\gamma_2}\log\left(\frac{\sigma_{\max}}{\sigma_0}\right),
\]
since $\gamma_2>1$. The desired result \req{unsucc-neg} then follows
from the equality $k +1 = |\calS_k| + |\calU_k|$
and the inequality $\gamma_1 < 1$ given by \req{eta-gamma2}.\hfill$\Box$\\

\noindent
{\bf Proof of Lemma~\ref{njsp-l-a}.}\label{derns-Lemma-proof}
We first observe that $\nabla_s^j \big(\|s\|^{p+\beta}\big)$ is a $j$-th order tensor,
whose norm is defined using \req{Tnorm}.  Moreover, using the relationships
\beqn{donce}
\nabla_s \big(\|s\|^\tau \big) = \tau\, \|s\|^{\tau-2}s
\tim{ and }
\nabla_s  \big(s^{\tau \otimes}\big) = \tau\, s^{(\tau-1)\otimes}\otimes I,
\ms (\tau \in \Re),
\eeqn
defining
\beqn{munu-def}
\nu_0 \eqdef 1,
\tim{ and }
\nu_i \eqdef \prod_{\ell=1}^{i}(p+2-2\ell),
\eeqn
and proceeding by induction, we obtain that, for some $\mu_{j,i}\geq 0$ with
$\mu_{1,1}=1$,
\[
\begin{array}{ll}
\nabla_s\left[\nabla_s^{j-1} \big(\| s\|^{p+\beta} \big) \right]&\\
& \hspace*{-3cm} = \nabla_s\left[ \bigsum_{i=2}^j \mu_{j-1,i-1} \nu_{i-1}
   \|s\|^{p+\beta-2(i-1)} \, s^{(2(i-1)-(j-1)) \otimes} \otimes I^{((j-1)-(i-1))\otimes} \right]\\
&\hspace*{-3cm} = \bigsum_{i=2}^j \mu_{j-1,i-1} \nu_{i-1} \Big[
  (p+\beta-2(i-1))\|s\|^{p+\beta-2(i-1)-2} \, s^{(2(i-1)-(j-1)+1) \otimes} \otimes I^{(j-i)\otimes}\\
&  \hspace*{-1cm} + ((2(i-1)-(j-1)) \|s\|^{p+\beta-2(i-1)} \, s^{(2(i-1)-(j-1)-1)\otimes}
  \otimes I^{(j-1)-(i-1)+1)\otimes} \Big]\\
&\hspace*{-3cm} = \bigsum_{i=2}^j \mu_{j-1,i-1} \nu_{i-1} \Big[
  (p+\beta+2-2i)\|s\|^{p+\beta-2i} \, s^{(2i-j) \otimes} \otimes I^{(j-i)\otimes}\\
  &  \hspace*{-1cm} + (2(i-1)-j+1) \|s\|^{p+\beta-2(i-1)} \, s^{(2(i-1)-j) \otimes}
  \otimes I^{(j-(i-1))\otimes} \Big]\\
&\hspace*{-3cm} = \bigsum_{i=2}^j \mu_{j-1,i-1} \nu_{i-1} 
  (p+\beta+2-2i)\|s\|^{p+\beta-2i} \, s^{(2i-j) \otimes} \otimes I^{(j-i)\otimes}\\
&  \hspace*{-1cm} + \bigsum_{i=1}^{j-1} (2i-j+1) \mu_{j-1,i}\nu_i \|s\|^{p+\beta-2i}\, s^{(2i-j) \otimes}
  \otimes I^{(j-i)\otimes} \Big]\\
& \hspace*{-3cm} = \bigsum_{i=1}^j\big((p+\beta+2-2i)\mu_{j-1,i-1}\nu_{i-1}
   + (2i-j+1)\mu_{j-1,i}\nu_i \big) \|s\|^{p+\beta-2i} \, s^{(2i-j) \otimes} \otimes I^{(j-i)\otimes}.  
\end{array}
\]
where the last equation uses the convention that $\mu_{j,0} = 0$ for all $j$.
Thus we may write  
\beqn{nablaj-ns-full}
\nabla_s^j \big(\|s\|^{p+\beta}\big)
=\nabla_s\left[\nabla_s^{j-1} \big(\| s\|^{p+\beta} \big) \right]
= \bigsum_{i=1}^j \mu_{j,i} \nu_i \,
   \|s\|^{p+\beta-2i} \, s^{(2i-j) \otimes} \otimes I^{(j-i)\otimes}
\eeqn
with
\beqn{rec}
\begin{array}{lcl}
\mu_{j,i}\nu_i
& = & (p+\beta+2-2i) \mu_{j-1,i-1}\nu_{i-1} + (2i-j+1) \mu_{j-1,i}\nu_i \\*[1.5ex]
& = & \big[\mu_{j-1,i-1} + (2i-j+1) \mu_{j-1,i}\big]\nu_i,
\end{array}
\eeqn
where we used the identity
\beqn{nuratio}
\nu_i = (p+\beta+2-2i)\nu_{i-1} \tim{ for } i = 1, \ldots, j
\eeqn
to deduce the second equality. Now \req{nablaj-ns-full} gives that
\[
\nabla_s^j \big(\|s\|^{p+\beta}\big)[v]^j
= \bigsum_{i=1}^j \mu_{j,i} \nu_i
\|s\|^{p+\beta-j} \, \left(\frac{s^Tv}{\|s\|}\right)^{2i-j} (v^Tv)^{j-i}.
\]
It is then easy to see that the maximum in \req{Tnorm} is achieved for
$v = s/\|s\|$, so that 
\beqn{c0}
\| \, \nabla_s^j \big(\| s\|^{p+\beta} \big) \,\|_{[j]}
=\left(\bigsum_{i=1}^j \mu_{j,i}  \nu_i \right) \|s\|^{p+\beta-j}
= \pi_j \|s\|^{p+\beta-j}.
\eeqn
with
\beqn{pimunu}
\pi_j \eqdef \bigsum_{i=1}^{j}\mu_{j,i}\,\nu_i.
\eeqn
Successively using this definition, \req{rec}, \req{nuratio}
(twice), the identity $\mu_{j-1,j} = 0$ and \req{pimunu} again, we then deduce that
\beqn{pi-ineqs}
\begin{array}{lcl}
\pi_j
& = & \bigsum_{i=1}^{j} \mu_{j-1,i-1}\nu_i + \bigsum_{i=1}^{j} (2i-j+1) \mu_{j-1,i}\big]\nu_i\\*[1.5ex]
& = & \bigsum_{i=1}^{j-1} \mu_{j-1,i}\nu_{i+1} + \bigsum_{i=1}^{j} (2i-j+1) \mu_{j-1,i}\big]\nu_i\\*[1.5ex]
& = & \bigsum_{i=1}^{j-1} \mu_{j-1,i}\big[ \nu_{i+1} + (2i-j+1) \nu_i\big]\\*[1.5ex]
& = & \bigsum_{i=1}^{j-1} \mu_{j-1,i}\big[ (p+\beta+2-2(i+1))\nu_i + (2i-j+1) \nu_i\big]\\*[1.5ex]
& = & (p+\beta+1-j) \bigsum_{i=1}^{j-1} \mu_{j-1,i}\,\nu_i\\*[1.5ex]
& = & (p+\beta+1-j) \pi_{j-1},
\end{array}
\eeqn
Since $\pi_1 = p+\beta$ from the first part of \req{donce}, we
obtain that $\pi_j = (p+\beta)!/(p-j+\beta)!$, which,
combined with \req{c0} and \req{pimunu}, gives \req{npsp-a}.
We obtain \req{npsp-b} from \req{c0} and \req{pimunu}, the observation that
$\pi_p= (p+\beta)!$ and \req{pi-ineqs} for $j=p+1$.\hfill$\Box$

\section{Proof of Lemmas in Section \ref{upperbound-s}}

\noindent
{\bf Proof of Lemma \ref{Dm-lemma}.}
(See \cite[Lemma~2.1]{BirgGardMartSantToin17})
Observe that, because of \req{descent2} and \req{model},
\[
0 \leq m_k(0) - m_k(s_k)
= T_p(x_k,0) - T_p(x_k,s_k) -\frac{\sigma_k}{p+1} \|s_k\|^{p+\beta}
\]
which implies the desired bound. Note that $s_k\neq 0$ as long as we can
satisfy condition \req{descent2}, and so \req{Dphi} implies \req{rhokdef2} is
well defined.\hfill$\Box$\\ 

\noindent
{\bf Proof of Lemma \ref{sigmaupper-lemma}.}
(See \cite[Lemma~2.2]{BirgGardMartSantToin17})
Assume that
\beqn{siglarge}
\sigma_k \geq \frac{L}{1-\eta_2}.
\eeqn
Using \req{resf} and \req{Dphi}, we may then deduce that
\[
|\rho_k - 1|
\leq \frac{|f(x_k+s_k) - T_p(x_k,s_k)|}{|T_p(x_k,0)-T_p(x_k,s_k)|}
\leq \frac{L}{\sigma_k}
\leq 1-\eta_2
\]
and thus that $\rho_k \geq \eta_2$. Then iteration $k$ is very successful in
that $\rho_k \geq \eta_2$ and $\sigma_{k+1}\leq \sigma_k$.  As a consequence,
the mechanism of the algorithm ensures that \req{sigmaupper} holds.\hfill$\Box$\\

\end{document}